\documentclass[12pt,a4paper,english,reqno]{amsart}
\usepackage{amsthm,amssymb,amstext}
\usepackage[T1]{fontenc}
\usepackage{lmodern,microtype,calrsfs}
\usepackage{mathtools,enumitem}
\usepackage[colorlinks=true,allcolors=blue]{hyperref}

\theoremstyle{plain}
\newtheorem{thm}{Theorem}[section]
\newtheorem{prp}[thm]{Proposition}
\newtheorem{lem}[thm]{Lemma}
\newtheorem{cor}[thm]{Corollary}
\theoremstyle{definition}

\newtheorem{dfn*}[thm]{Definition}
\newtheorem{exm}[thm]{Example}
\newtheorem{rmk}[thm]{Remark}

\numberwithin{equation}{section}

\newcommand{\thmref}[1]{Theorem~\ref{#1}}
\newcommand{\prpref}[1]{Proposition~\ref{#1}}
\newcommand{\corref}[1]{Corollary~\ref{#1}}
\newcommand{\lemref}[1]{Lemma~\ref{#1}}
\def\al{\alpha}
\def\be{\beta}
\def\ga{\gamma}
\def\de{\delta}
\def\ep{\varepsilon}
\def\et{\eta}
\def\io{\iota}
\def\la{\lambda}
\def\ro{\rho}
\def\si{\sigma}
\def\ta{\tau}
\def\Ga{\Gamma}
\def\De{\Delta}
\def\La{\Lambda}
\def\Ph{\Phi}
\def\Th{\Theta}
\def\Hom{\operatorname{Hom}}
\def\rad{\operatorname{rad}}
\def\End{\operatorname{End}}
\def\prj{\operatorname{prj}}
\def\thick{\operatorname{thick}}
\def\add{\operatorname{add}}
\def\Kb{{\mathcal K}^{\text{\rm b}}}
\def\Db{{\mathcal D}^{\text{\rm b}}}
\def\perf{\operatorname{per}}
\def\calC{{\mathcal C}}
\def\calD{{\mathcal D}}
\def\bbZ{{\mathbb Z}}
\def\bbC{{\mathbb C}}
\def\incl{\hookrightarrow}
\def\iso{\cong}
\def\ds{\oplus}
\def\ox{\otimes}
\def\Lox{\overset{\mathbf{L}}{\otimes}}
\def\Ds{\bigoplus}
\def\ya#1{\xrightarrow{#1}}
\newcommand{\Fgt}{\operatorname{Fgt}}
\def\RHom{\operatorname{RHom}}
\def\REnd{\operatorname{REnd}}
\def\Jac{\operatorname{Jac}}
\def\gr{\operatorname{gr}}
\def\Cov{\operatorname{Cov}}
\def\Cone{\operatorname{Cone}}
\def\Loc{\operatorname{Loc}}
\newcommand{\CT}{\operatorname{CT}}

\def\Ainf{A_\infty}
\def\fd{\mathrm{fd}}
\newcommand{\uRHom}{\underline{\RHom}}
\newcommand{\uREnd}{\underline{\REnd}}

\begin{document}
\title[Compatible gradings and cluster equivalence]{Compatible gradings on dg models and cluster equivalence}
\author{Shengyong Pan}
\address{School of Mathematics and Statistics, Beijing Jiaotong University,
Beijing, 100044, China.}
\address{Beijing Key Laboratory of Biological Big Data and Topological Statistics, Beijing Jiaotong University,
Beijing, 100044, China.}
\email{shypan@bjtu.edu.cn}

\date{31 October 2026}
\begin{abstract}
Let $\La_1$ and $\La_2$ be finite-dimensional $\ta_2$-finite complex algebras of global dimension at most two. We prove that, if they are cluster equivalent, then they admit $\bbZ$-gradings for which the bounded derived categories of graded modules are triangle equivalent. The equivalence is not required to commute with the internal grading shifts. The proof uses common multigraded dg models, relative gradings on two-term silting objects, and covering quotients. We also consider tilting lifts of cluster-tilting objects with endomorphism algebras of global dimension at most two. In acyclic cluster type, the existence of such a lift for every cluster-tilting object implies piecewise heredity. We construct a counterexample without this assumption using a surface algebra associated with a one-holed torus.
\end{abstract}
\subjclass[2020]{16E35, 16G20, 18G80}
\keywords{Cluster category, graded derived equivalence, silting object,
dg algebra, torus action, Calabi--Yau completion, surface algebra}
\maketitle
\tableofcontents

\section{Introduction}
Derived equivalences of finite-dimensional algebras of global dimension at most two induce equivalences of their generalized cluster categories. The converse does not hold in general. Amiot and Oppermann studied this problem by introducing gradings and comparing derived categories of graded modules. The purpose of this paper is to study the existence of compatible gradings and the lifting of cluster-tilting objects to tilting complexes.

Throughout the paper, unless otherwise stated, we work over $\bbC$ and use right modules. Let $\La$ be a finite-dimensional algebra of global dimension at most two, and let $\mathbb S$ be the Serre functor of $\Db(\La)$. Put $\mathbb S_2=\mathbb S[-2]$. Amiot's generalized cluster category is the triangulated hull of the orbit category $\Db(\La)/\mathbb S_2$; it admits the dg description
\[
 \calC_\La
 =\perf\Pi_3(\La)/\calD_{\fd}(\Pi_3(\La)),
\]
where $\Pi_3(\La)$ is the derived $3$-Calabi--Yau completion; see \cite{Amiot,Keller}. We assume that $\La$ is $\ta_2$-finite, that is, $H^0\Pi_3(\La)$ is finite dimensional. Under this hypothesis, $\calC_\La$ is Hom-finite and $2$-Calabi--Yau, and the canonical projection $\pi:\Db(\La)\to\calC_\La$ sends $\La$ to a cluster-tilting object. Two such algebras are called \emph{cluster equivalent} when their generalized cluster categories are triangle equivalent.

The endomorphism algebra of the canonical cluster-tilting object has a natural grading
\[
 \widetilde\La
 =\End_{\calC_\La}(\pi\La)
 \iso\Ds_{n\geq0}
       \Hom_{\Db(\La)}(\La,\mathbb S_2^{-n}\La),
 \qquad (\widetilde\La)_0\iso\La.
\]
Thus $\La$ is the degree-zero subalgebra of $\widetilde\La$. We use gradings of this kind to compare covering categories of cluster equivalent algebras.

Under a compatibility assumption, Amiot and Oppermann obtained an equivalence of graded derived categories \cite[Theorem~8.7]{AO}. This leads to the following question: do any two cluster equivalent algebras admit $\bbZ$-gradings for which their bounded derived categories of graded modules are triangle equivalent? The gradings in this question are not fixed.

A given cluster equivalence need not preserve the canonical cluster-tilting objects. We therefore choose a silting lift of the image of one of these objects and consider its endomorphism dg algebra. The comparison obtained in this way is ungraded. To retain both gradings, we construct commuting torus actions on a common minimal $\Ainf$ model.

Our main result is the following.

\begin{thm}\label{thm:main}
Let $\La_1$ and $\La_2$ be finite-dimensional $\ta_2$-finite $\bbC$-algebras of global dimension at most two. If $\calC_{\La_1}$ and $\calC_{\La_2}$ are $\bbC$-linear triangle equivalent, then there are $\bbZ$-gradings on $\La_1$ and $\La_2$ for which
\[
 \Db(\gr\La_1)\simeq\Db(\gr\La_2)
\]
as $\bbC$-linear triangulated categories.
\end{thm}

The equivalence in \thmref{thm:main} is an equivalence of triangulated categories; compatibility with the internal grading shifts is not assumed. Nor is it required to induce a given cluster equivalence. Under shift compatibility, the underlying algebras are derived equivalent (\prpref{prop:shift-descent}). Without it, the converse of \thmref{thm:main} fails (Example~\ref{ex:no-converse}).

For the proof of \thmref{thm:main}, we establish the following results. In \thmref{thm:common-dg}, we construct a common bigraded dg model. In \thmref{thm:relative}, we prove that a two-term silting object admits an additional grading without changing its given graded realization. The resulting models are compared by \prpref{prop:cover} and \thmref{thm:automatic}. The ordinary dg identification follows from \cite[Corollary~4.27]{HK} and noncommutative Mather--Yau \cite[Theorems~3.5 and~4.3]{HZ}, using the vanishing of the canonical potential class. This identification is used without a grading. The distinction from graded recognition is discussed after \prpref{prop:zero-class-recognition}; see \cite[Example~7.6]{HZFinite}.

An algebra has \emph{acyclic cluster type} if its cluster category is equivalent to $\calC_{\bbC Q}$ for a finite acyclic quiver $Q$. This refers to the cluster category, rather than to the ordinary quiver of the algebra. Amiot and Oppermann described derived equivalence in this setting using graded mutation \cite[Theorem~3.15]{AOTree}. An algebra of tree cluster type is derived equivalent to a hereditary path algebra (\cite[Corollary~3.16]{AOTree}), whereas in type $\widetilde A$ the derived classes are distinguished by an integer weight, up to sign in the symmetric case (\cite[Theorems~4.5 and~4.7]{AOTree}). These results also motivate the second theme of this paper.

We also study the lifting of cluster-tilting objects to tilting complexes whose endomorphism algebras have global dimension at most two. The second part of \cite[Question 6.23, p.~47]{AOTree} asks whether lifting every cluster-tilting object forces the original algebra to be piecewise hereditary. The question itself has no acyclic-cluster-type assumption. We prove the implication with this additional hypothesis by \prpref{prop:acyclic-lifting}: it suffices to lift one object whose endomorphism algebra is hereditary. In \thmref{thm:torus-counterexample}, we construct a counterexample in the general setting. Its algebra has an acyclic ordinary quiver but does not have acyclic cluster type. The surface argument is independent of \thmref{thm:main}, the common dg model and noncommutative Mather--Yau.

The paper is organized as follows. In Section~\ref{sec:background}, we fix the notation and recall the dg and cluster constructions. In Section~\ref{sec:relative}, we add a grading to a two-term silting object while retaining its first graded realization. In Section~\ref{sec:common}, we construct the common multigraded dg model. We supply the ordinary dg comparison needed for its application in Section~\ref{sec:completion}, and in Section~\ref{sec:cover} we identify the covering quotients. These results are combined in Section~\ref{sec:main-proof} to prove \thmref{thm:main}. In Section~\ref{sec:scope} we treat internal shifts and the converse question. Finally, in Section~\ref{sec:torus-lifting} we consider tilting lifts in acyclic cluster type and on the marked torus.

\section*{Acknowledgments}
The mathematical ideas in this paper are the author's own.
The author used ChatGPT, to help draft and revise the manuscript. The author is responsible
for the mathematical content of the paper.
This work is supported by the Beijing Natural Science Foundation (Grant Nos. 1262017, 1252011).

\section{Preliminaries}\label{sec:background}
Throughout the paper, $k=\bbC$ unless otherwise stated, and all modules
are right modules. We fix the notation for dg modules and internal
gradings, recall the constructions used in the proof, and collect the
required comparison results. We use the notation for graded categories and coverings of \cite{AP},
with the internal shift convention specified below.

For a class $\mathcal S$ in an additive category, $\add\mathcal S$
denotes its closure under finite direct sums and direct summands.
In a triangulated category, $\thick(\mathcal S)$ also allows shifts
and cones, whereas $\Loc(\mathcal S)$ allows arbitrary coproducts.
For an ordinary algebra $L$, we write $\operatorname{mod}L$ for
finitely generated right modules, $\prj L$ for finitely generated
projective right modules, and $\Kb(\prj L)$ for their bounded homotopy
category. Complexes are cohomological:
\[
 X[n]^i=X^{i+n},\qquad d_{X[n]}=(-1)^n d_X.
\]
Only cohomological degrees enter the Koszul sign rule.
We write $D=\Hom_k(-,k)$ for the $k$-linear dual and
$\Db(L)=\Db(\operatorname{mod}L)$. Unless a base algebra is indicated,
tensor products are taken over $k$.

\subsection{Dg categories, modules and Morita equivalences}\label{subsec:dg-definitions}

\begin{dfn*}[{\cite[Sections 2.2--2.3]{KellerDG}}]\label{def:dg-categories}
A dg $k$-category $\mathcal A$ has complexes $\mathcal A(X,Y)$,
associative composition maps of complexes
\[
 \mathcal A(Y,Z)\ox_k\mathcal A(X,Y)\longrightarrow\mathcal A(X,Z),
 \qquad f\ox g\longmapsto fg,
\]
and closed degree-zero identities $1_X$. In particular,
\[
 d(fg)=d(f)g+(-1)^{|f|}f\,d(g).
\]
A dg functor $F:\mathcal A\to\mathcal B$ consists of a map on
objects and maps of complexes
$F_{X,Y}:\mathcal A(X,Y)\to\mathcal B(FX,FY)$ satisfying
$F(fg)=F(f)F(g)$ and $F(1_X)=1_{FX}$.
\end{dfn*}

A dg algebra $R$ is a dg category with one object $*$ and
$\mathcal A(*,*)=R$. Thus $R=\Ds_iR^i$ is a unital graded algebra
with $d_R^2=0$ and the displayed Leibniz rule. We write
\[
 \operatorname{Ob}(H^0\mathcal A)=\operatorname{Ob}(\mathcal A),
 \qquad (H^0\mathcal A)(X,Y)=H^0\mathcal A(X,Y).
\]

\begin{dfn*}[{\cite[Section 2.3]{KellerDG}}]\label{def:quasi-equivalence}
A dg functor $F:\mathcal A\to\mathcal B$ is a
\emph{quasi-equivalence} if every $F_{X,Y}$ is a quasi-isomorphism
and $H^0F$ is essentially surjective.
\end{dfn*}

Let $\mathcal C_{\mathrm{dg}}(k)$ be the dg category of complexes,
with
\[
 \mathcal C_{\mathrm{dg}}(k)(X,Y)^n
 =\prod_i\Hom_k(X^i,Y^{i+n}),\qquad
 d(f)=d_Yf-(-1)^nf d_X.
\]
The opposite dg category has
$\mathcal A^{\mathrm{op}}(X,Y)=\mathcal A(Y,X)$ and composition
$f\circ_{\mathrm{op}}g=(-1)^{|f||g|}gf$.

\begin{dfn*}[{\cite[Section 3.1]{KellerDG}}]\label{def:dg-modules}
A right dg module over a small dg category $\mathcal A$ is a dg
functor
\[
 M:\mathcal A^{\mathrm{op}}\longrightarrow\mathcal C_{\mathrm{dg}}(k).
\]
Morphisms of dg modules are closed degree-zero dg natural
transformations.
\end{dfn*}

For a dg algebra $R$, this is the usual right dg module: a graded
right $R$-module $X$ with
\[
 d_X^2=0,\qquad d_X(xa)=d_X(x)a+(-1)^{|x|}x\,d_R(a).
\]
We use the right-module convention
\[
 \begin{gathered}
\Hom_R^n(X,Y)=\{f:X\to Y:|f|=n,\ f(xa)=f(x)a\},\\
 d(f)=d_Yf-(-1)^nf d_X.
\end{gathered}
\]
Thus $K(R)$ has morphism spaces $H^0\Hom_R^\bullet(X,Y)$, and
$\calD R=K(R)[\mathrm{qis}^{-1}]$; see
\cite[Section 3.2]{KellerDG}. We put
\[
 \perf R=\thick(R),\qquad
 \calD_{\fd}R=\{X\in\calD R:\sum_i\dim_kH^iX<\infty\}.
\]
The notation $R\simeq S$ means that $R$ and $S$ are connected
by a zigzag of quasi-isomorphisms of dg algebras, unless an actual map is specified.

A dg module $P$ is homotopically projective if
$\Hom_R^\bullet(P,Z)$ is acyclic whenever $Z$ is acyclic.
Semi-free resolutions, whose filtration quotients are sums of shifts
of $R$, provide such models. We use
$\REnd_R(X)=\End_R^\bullet(P)$ for a semi-free resolution $P\to X$,
with multiplication given by composition.

\begin{dfn*}[{\cite[Section 3.8]{KellerDG}}]\label{def:dg-morita}
An $A$--$R$ dg bimodule $P$ induces a dg Morita equivalence if
\[
 -\Lox_AP:\calD A\rightleftarrows\calD R:\RHom_R(P,-)
\]
are mutually inverse equivalences. The bimodule and its two actions
are part of the data.
\end{dfn*}

\begin{dfn*}[{\cite[Section 2]{Amiot}; \cite[Section 3.3]{Keller}}]\label{def:smooth-cy}
Put $R^e=R^{\mathrm{op}}\ox_kR$. The dg algebra $R$ is
homologically smooth if $R\in\perf(R^e)$. It is bimodule
$3$-Calabi--Yau if it is smooth and
\[
 \RHom_{R^e}(R,R^e)\simeq R[-3]
\]
as dg bimodules, using the standard interchange of the two factors
to view the dual as a bimodule.
\end{dfn*}

We call $R$ connective if $H^{>0}R=0$, nonpositive if $R^{>0}=0$,
and cohomologically locally finite if each $H^iR$ is finite-dimensional.
These terms refer to cohomological degrees, after forgetting internal
gradings. A pseudocompact algebra is a complete Hausdorff topological
algebra with a basis of zero consisting of ideals of finite codimension;
equivalently, it is an inverse limit of finite-dimensional discrete
algebras. For dg algebras we use the degreewise version and continuous
maps, as in \cite[Section 6 and Appendix A.2]{VdB}.
We denote forgetting the topology by $\Fgt_{\mathrm{top}}$;
this does not change the cohomological grading or differential.

\subsection{Internal gradings and torus actions}\label{subsec:gradings}

\begin{dfn*}[{\cite[Sections 4.1--4.2]{KalckYang}}]\label{def:internal-gradings}
Let $G=\bbZ^r$. An internally $G$-graded dg algebra is a dg algebra
with a direct-sum decomposition $R=\Ds_{i,g}R^i_g$ such that
\[
 1\in R^0_0,\qquad R^i_gR^j_h\subset R^{i+j}_{g+h},\qquad
 d_R(R^i_g)\subset R^{i+1}_g.
\]
A graded dg module is a dg module $X=\Ds_{i,g}X^i_g$ satisfying
$X^i_gR^j_h\subset X^{i+j}_{g+h}$ and
$d_X(X^i_g)\subset X^{i+1}_g$.
A morphism of graded dg modules preserves the internal degree.
\end{dfn*}

For $G=\bbZ^2$, we write $R^i_{p,q}$ and call $R$ a bigraded dg
algebra. Thus there are three degrees: one cohomological and two
internal (Adams) degrees. The differential has degree $(1,0,0)$.
More generally, for $G=\bbZ^r\ds\bbZ^s$, we write
$X^i_{p,q}$, where $p\in\bbZ^r$ and $q\in\bbZ^s$.
For $g,h\in G$, the internal shift is defined by
$X(g)^i_h=X^i_{h+g}$. In two grading blocks this convention becomes
\[
 X(a,b)^i_{p,q}=X^i_{p+a,q+b},\qquad
 (\Fgt_qR)^i_p=\Ds_qR^i_{p,q}.
\]
The functor $\Fgt_q$ forgets the $q$-grading on algebras, modules
and morphisms; $\Fgt_p$ is defined similarly. We write $\Fgt$ when
all internal gradings are forgotten. These functors leave the
differential unchanged. Forgetting a grading is
not the totalization of \cite[Section 4.1]{KalckYang}.

Following \cite[Sections 4.1--4.2]{KalckYang}, let $\calD^G R$
be the localization of the graded homotopy category at graded
quasi-isomorphisms, and set
\[
 \begin{aligned}
 \perf^G R&=\thick\{R(g):g\in G\},\\
 \calD^G_{\fd}R&=\{X:\sum_{i,g}\dim_kH^i_gX<\infty\}.
 \end{aligned}
\]
Here $H^i_gX$ denotes the internal degree-$g$ component of $H^iX$.
If $R$ has grading blocks $p\in\bbZ^r$ and $q\in\bbZ^s$, we use
\[
 \calD^{\bbZ^r}_pR:=\calD^{\bbZ^r}(\Fgt_qR),\qquad
 \perf^{\bbZ^r}_pR:=\perf^{\bbZ^r}(\Fgt_qR).
\]
Thus the subscript $p$ specifies the grading retained, and the
superscript specifies its grading group. We also write
$\calD^pR$ and $\perf^pR$ for these categories when the rank is clear.
The notation with $q$ in place of $p$ is defined in the same way.
If $R$ has only a $p$-grading, no grading is forgotten.

For an ordinary graded algebra $L$, $\gr L$ denotes finitely generated
graded modules and degree-zero maps; $\operatorname{Gr}L$ allows all
graded modules. We write $\gr^pL$ or $\gr^qL$ when we need to
specify the grading used. For example, $\gr^qL$ consists of finitely
generated modules $M=\Ds_qM_q$ with $M_qL_b\subseteq M_{q+b}$,
and morphisms preserving $q$. The superscript denotes a grading,
not a homogeneous component. We write $\gr^p\prj L$ for its full
subcategory of finitely generated graded projective modules and
$\Kb(\gr^p\prj L)$ for their bounded homotopy category.
For a grading group $G$, we also use $\gr_G\prj L$.

The homogeneous Hom complex is defined by
\[
 \underline{\Hom}_R(X,Y)=
 \Ds_{g\in G}\Hom_{R,\mathrm{gr}}^\bullet(X,Y(g)).
\]
An element in its cohomological degree $n$ and internal degree $g$
sends $X^i_h$ into $Y^{i+n}_{h+g}$. For a graded semi-free resolution
$P\to X$, we use the models
\[
 \uRHom_R(X,Y)=\underline{\Hom}_R(P,Y),\qquad
 \uREnd_R(X)=\underline{\End}_R(P).
\]
Here the multiplication is induced by composition of homogeneous
endomorphisms. To specify the grading group, we write
$\uRHom_R^G(X,Y)$ and $\uREnd_R^G(X)$. In particular,
\[
 \bigl(\uREnd_R^{\bbZ^2}(X)\bigr)^n_{a,b}
 =\Hom^n_{R,\mathrm{gr}}(P,P(a,b)).
\]
Thus the underline means that all homogeneous internal degrees are
included; it does not mean only endomorphisms of internal degree zero.
For the $p$-grading alone, we write $\uREnd_R^p(X)$, with $R$
interpreted as $\Fgt_qR$ if it originally has both gradings.
The same conventions apply to $\uRHom$ and to grading blocks of
arbitrary finite rank.

For a bigraded dg algebra $R$, the notation $(R,p)$ means
$\Fgt_qR$, and $(R,q)$ means $\Fgt_pR$. We write
\[
 (R_{p=0})^i_q=R^i_{0,q},\qquad
 (R_{q=0})^i_p=R^i_{p,0}.
\]
These are dg subalgebras with the remaining internal grading.
When the completion grading is denoted by $u$, we abbreviate
$R_{u=0}$ to $R_0$. This notation is distinct from $R^0$, the
cohomological degree-zero component, and from $H^0R$.

For graded algebras $L_1,L_2$, a triangle equivalence
$F:\Db(\gr L_1)\to\Db(\gr L_2)$ is not required to commute
with the internal shifts. We say that $F$ is
\emph{compatible with the grading shifts} when it is equipped
with natural isomorphisms
\[
 \eta_a:F\circ(a)\xrightarrow{\sim}(a)\circ F
 \qquad(a\in\bbZ)
\]
satisfying $\eta_0=1$ and
$\eta_{a+b,X}=(a)(\eta_{b,X})\circ\eta_{a,X(b)}$.
For a triangle equivalence, compatibility with the cohomological
shift $[1]$ is part of its triangle functor structure.

\begin{dfn*}[{\cite[Expos\'e I, Section 4]{SGA}}]
\label{def:torus-action}
Let $k$ be a field and let $r\geq 0$. The split torus of rank $r$
over $k$ is the affine group scheme
\[
 T=(\mathbb G_{m,k})^r.
\]
Its coordinate Hopf algebra is
$k[T]=k[t_1^{\pm1},\ldots,t_r^{\pm1}]$,
with comultiplication $\De$ and counit $\epsilon$ defined by
\[
 \De(t_j)=t_j\ox t_j,
 \qquad
 \epsilon(t_j)=1
 \qquad (1\leq j\leq r).
\]
All tensor products in this definition are taken over $k$.

Let $V$ be a $k$-vector space, not necessarily finite dimensional.
A \emph{rational $T$-action} on $V$ is a $k$-linear map
\[
 \de:V\longrightarrow V\ox k[T]
\]
such that
\[
 (\de\ox 1_{k[T]})\de=(1_V\ox\De)\de,
 \qquad
 (1_V\ox\epsilon)\de=1_V,
\]
where $V\ox k$ is identified with $V$.
The pair $(V,\de)$ is called a \emph{rational $T$-representation}.

Let $R$ be a unital dg $k$-algebra. Regard $k[T]$ as a dg algebra
concentrated in cohomological degree zero with zero differential.
A \emph{rational $T$-action} on $R$ is a rational $T$-action on
its underlying vector space for which
\[
 \de:R\longrightarrow R\ox k[T]
\]
is a morphism of unital dg $k$-algebras. Thus $\de$ preserves
cohomological degrees and satisfies
\[
 \de(1_R)=1_R\ox1,\qquad
 \de(ab)=\de(a)\de(b),\qquad
 \de d_R=(d_R\ox1)\de.
\]
\end{dfn*}

The character lattice $X^*(T)=\Hom_{k\text{-group schemes}}(T,\mathbb G_{m,k})$
is identified with $\bbZ^r$. For $g=(g_1,\ldots,g_r)$, put
$t^g=t_1^{g_1}\cdots t_r^{g_r}$. Rationality implies
\[
 V=\Ds_{g\in\bbZ^r}V_g,\qquad
 V_g=\{v:\de(v)=v\ox t^g\}.
\]
Thus rational torus actions on dg algebras are equivalent to internal
gradings. Equivariant maps preserve weights, and taking a fixed weight
is exact. For two commuting torus actions, we use the weight decomposition
for their product action. Over $k=\bbC$ this reads
$\ro(s,t)x=s^pt^qx$ for $x\in R^i_{p,q}$.
For a completed algebra, a continuous action rational on a cofinal
system of finite-dimensional quotients need not make the whole algebra
a direct sum of weight spaces. This distinction is used in
Subsection~\ref{subsec:localfinite}.

\subsection{Minimal models}\label{subsec:dg-properties}

\begin{dfn*}[{\cite[Sections 3.1 and 3.4--3.5]{KellerAinf}}]\label{def:ainfinity}
An uncurved $\Ainf$ algebra is a graded space $H$ with maps
$m_n:H^{\ox n}\to H$ of degree $2-n$, for $n\ge1$, satisfying
\[
 \sum_{r+s+t=n}(-1)^{r+st}
 m_{r+1+t}(1^{\ox r}\ox m_s\ox1^{\ox t})=0
 \quad(r,t\ge0,\ s\ge1),
\]
with the Koszul convention for tensor products of maps.
It is minimal if $m_1=0$. It is strictly unital if there is
$1\in H^0$ such that $m_2(1,x)=x=m_2(x,1)$ and
$m_n(\ldots,1,\ldots)=0$ for $n\ne2$.
\end{dfn*}

We use the equivalent coalgebra description from
\cite[Section 3.6]{KellerAinf}. Put $V=H[1]$ and
\[
 T^c(V)=\Ds_{n\ge0}V^{\ox n},\qquad
 \De(v_1\cdots v_n)=\sum_{j=0}^n
       (v_1\cdots v_j)\ox(v_{j+1}\cdots v_n).
\]
The suspended operations determine a degree-one coderivation $b$
with $b^2=0$ and zero constant term. An $\Ainf$ morphism is a
degree-zero coaugmented coalgebra map commuting with coderivations.
Its Taylor components $F_n$ have unsuspended degrees $1-n$.
It is a quasi-isomorphism if $F_1$ is a quasi-isomorphism, and is
strictly unital if $F_1(1)=1$ and the higher components vanish on
tensors containing $1$.

\begin{dfn*}[{\cite[Section 3.3]{KellerAinf}}]
\label{def:minimal-model}
Let $k$ be a field. An $\Ainf$-algebra $(H,(m_n)_{n\geq1})$
is called \emph{minimal} if $m_1=0$.

Let $R$ be a dg $k$-algebra. A \emph{minimal model} of $R$
consists of an $\Ainf$-algebra structure $(m_n)_{n\geq1}$
on $H^*R$ and an $\Ainf$ quasi-isomorphism
\[
 f=(f_n)_{n\geq1}:(H^*R,(m_n))\longrightarrow R
\]
such that
\[
 m_1=0,\qquad
 m_2([a],[b])=[ab],\qquad
 H^*(f_1)=1_{H^*R}.
\]
Here $a,b$ are cocycles in $R$, and $R$ is regarded as an
$\Ainf$-algebra with $m_1=d_R$, $m_2$ its multiplication,
and $m_n=0$ for $n\geq3$.

The dg algebra $R$ is called \emph{formal} if it admits
a minimal model for which
\[
 m_n=0\qquad(n\geq3).
\]
\end{dfn*}

By Kadeishvili's theorem
\cite[Section 3.3]{KellerAinf}, every dg algebra over a field
admits a minimal model, unique up to a not necessarily unique
$\Ainf$-isomorphism.

Suppose that $R$ has an internal grading by an abelian group $G$.
A minimal model $(H^*R,(m_n),f)$ is called \emph{internally graded}
if, for every $n\geq1$ and $g_1,\ldots,g_n\in G$,
\[
\begin{aligned}
 m_n\bigl((H^*R)_{g_1}\ox\cdots\ox(H^*R)_{g_n}\bigr)
 &\subseteq (H^*R)_{g_1+\cdots+g_n},\\
 f_n\bigl((H^*R)_{g_1}\ox\cdots\ox(H^*R)_{g_n}\bigr)
 &\subseteq R_{g_1+\cdots+g_n}.
\end{aligned}
\]
Thus all structure maps and all components of $f$ have
internal degree zero.

\subsection{Properties of dg models}\label{subsec:standard-properties}
We recall the following properties.

\begin{enumerate}[label=(\roman*)]
\item A compact object is one whose Hom functor commutes with
coproducts. The compact objects of $\calD R$ and $\calD^G R$ are,
respectively, $\perf R$ and $\perf^G R$
\cite[Theorem 3.4 and Corollary 3.7]{KellerDG}. A graded dg
quasi-isomorphism induces inverse extension and restriction
equivalences, preserving regular objects, internal shifts and
finite-dimensional total cohomology \cite[Section 3.8]{KellerDG}.

\item The Morita criterion \cite[Lemma 3.10]{KellerDG} says that a
bimodule $P$ as in Definition~\ref{def:dg-morita} induces an equivalence
if $P\in\perf R$, $\thick(P)=\perf R$, and its action map
$A\to\REnd_R(P)$ is a quasi-isomorphism. In the graded case the
conditions are
\[
 P\in\perf^G R,\qquad
 \thick\{P(g):g\in G\}=\perf^G R,\qquad
 A\ya{\sim}\uREnd_R(P).
\]
The resulting graded equivalence commutes with internal shifts.

\item If $L$ is finite-dimensional of finite global dimension, then
\[
 \perf L\simeq\Db(\operatorname{mod}L),\qquad
 \perf^{\bbZ}L\simeq\Db(\gr L).
\]
For the graded assertion, take finite graded free resolutions; the
degree-zero component of an ungraded splitting splits a graded
projective syzygy. If $R$ is smooth, then
$\calD_{\fd}R\subset\perf R$: tensor a finite-cell model of the
diagonal, up to retract, with $X\in\calD_{\fd}R$. Its cells become
$X\ox_kR$, which are perfect. Compare
\cite[Section 2, before Theorem 2.1]{Amiot}.

\item A connective dg algebra has the nonpositive model
\[
 (\ta_{\le0}R)^i=
 \begin{cases}R^i&i<0,\\Z^0R&i=0,\\0&i>0,\end{cases}
 \qquad \ta_{\le0}R\ya{\sim}R.
\]
This construction preserves internal gradings.

\item Minimal models exist over a field and are unique up to
$\Ainf$ isomorphism \cite[Sections 3.3 and 8.1]{KellerAinf}.
Choose decompositions $R^i=B^i\ds\mathcal H^i\ds L^i$, where
$Z^iR=B^i\ds\mathcal H^i$ and $d:L^i\iso B^{i+1}$.
We thus obtain contraction data
\[
 \pi\io=1,\qquad 1-\io\pi=dh+hd,\qquad
 h^2=h\io=\pi h=0.
\]
The transfer formulas use multiplication at vertices, $h$ on internal
edges, $\io$ at leaves and $\pi$ at the root of planar trees.
Homogeneous choices make every formula homogeneous. If the cohomology
unit is nonzero, take $1\in\mathcal H^0$ to retain strict units.

\item The unital enveloping dg algebra $U(H)$ satisfies
\[
 H\ya{\sim}U(H),\qquad
 \Hom_{\mathrm{dgAlg}}(U(H),R)\iso
 \Hom_{\Ainf,\mathrm{su}}(H,R)
\]
by \cite[Proposition 2.1]{COS}. The first map is an $\Ainf$
quasi-isomorphism. The construction needs no augmentation and retains
internal gradings. If $H^{>0}=0$ and each $H^i$ is finite-dimensional,
then, for $V=H[1]$,
\[
 C_N=\Ds_{-N\le j\le0}T^c(V)^j,
 \qquad \dim C_N<\infty,\qquad T^c(V)=\bigcup_N C_N.
\]
These are dg subcoalgebras: $V$ is strictly negative, so tensor length
is at most $N$ in $C_N$. Compatibility of automorphisms across these
subcoalgebras is the separate assertion of Proposition~\ref{prop:torus}.
\end{enumerate}

\subsection{Orbit and generalized cluster categories}\label{subsec:cluster-definitions}

\begin{dfn*}[{\cite[Section 1]{KellerOrbit}}]\label{def:orbit-category}
For an automorphism $F$ of a $k$-linear category $\mathcal T$, the
orbit category $\mathcal T/F$ has the same objects as $\mathcal T$,
with
\[
 \Hom_{\mathcal T/F}(X,Y)=\Ds_{n\in\bbZ}\Hom_{\mathcal T}(X,F^nY).
\]
For $f_n:X\to F^nY$ and $g_m:Y\to F^mZ$, composition is
$F^n(g_m)f_n$.
\end{dfn*}

An autoequivalence can be replaced by an equivalent strict action.
The orbit category need not be triangulated. From a specified dg enhancement and a lift of $F$, one constructs
a dg orbit category and its triangulated hull;
see \cite[Sections 4--5]{KellerOrbit} and \cite{KellerOrbitCorr}.

\begin{dfn*}[{\cite[Section 4.1]{Amiot}}]\label{def:generalized-cluster}
Let $\La$ be a finite-dimensional algebra of finite global dimension,
and put $B_\La=\La\ds D\La[-3]$, where $D=\Hom_k(-,k)$.
Give $B_\La$ zero differential and the trivial-extension product
$(a,f)(b,g)=(ab,ag+fb)$. If $p^*$ is restriction along
$B_\La\twoheadrightarrow\La$, set
\[
 \mathcal E_\La=
 \thick_{\calD_{\fd}B_\La}(p^*\Db(\operatorname{mod}\La)),
 \qquad \calC_\La=\mathcal E_\La/\perf B_\La.
\]
The category $\calC_\La$ is the generalized cluster category of $\La$.
\end{dfn*}

The sequence $0\to D\La[-3]\to B_\La\to\La\to0$ ensures that
the denominator lies in $\mathcal E_\La$. The canonical functor is
$p^*$ followed by the quotient. This realizes the triangulated hull
of $\Db(\operatorname{mod}\La)/\mathbb S[-2]$, where
$\mathbb S=-\Lox_\La D\La$; one obtains the
same orbit by using $\mathbb S^{-1}[2]$. The image-generated numerator is essential; see
\cite[Section 1]{KellerOrbitCorr}.

\begin{dfn*}[{\cite[Section 4.1]{Keller}}]\label{def:cy-completion}
For a homologically smooth algebra $\La$, choose a cofibrant bimodule
$\Th\simeq\RHom_{\La^e}(\La,\La^e)[2]$.
Its undeformed derived $3$-Calabi--Yau completion is
\[
 \Pi_3(\La)=T_\La\Th=\Ds_{n\ge0}\Th^{\ox_\La n}.
\]
The canonical internal grading assigns degree $n$ to the $n$th tensor
power. The tensor algebra here is an ordinary direct sum.
\end{dfn*}

For finite-dimensional $\La$ of global dimension at most two,
write $\Ga=\Pi_3(\La)$. Then
$H^0\Ga\iso T_\La\operatorname{Ext}^2_\La(D\La,\La)$.
The algebra $\La$ is called $\ta_2$-finite when this algebra is
finite-dimensional \cite[Section 2]{AO}. Under this hypothesis,
\[
 \calC_\La\simeq\perf\Ga/\calD_{\fd}\Ga,
 \qquad \End_{\calC_\La}(\pi\Ga)\iso H^0\Ga.
\]
The category is Hom-finite and $2$-Calabi--Yau, meaning that there are
bifunctorial isomorphisms
\[
 D\Hom(X,Y)\iso\Hom(Y,X[2]).
\]
These statements,
including the canonical cluster-tilting object $\pi\Ga$, are
\cite[Theorem 2.1, Proposition 4.9 and Theorem 4.10]{Amiot}.

\begin{dfn*}[{\cite[Section 2]{Amiot}}]\label{def:cluster-tilting}
In a Hom-finite $2$-Calabi--Yau triangulated category, an object $T$
is cluster-tilting if $\add T$ is functorially finite and
\[
 \add T=\{X:\Hom(T,X[1])=0\}
        =\{X:\Hom(X,T[1])=0\}.
\]
Functorial finiteness means that every object admits both a right and
a left $\add T$-approximation. A right approximation $T_X\to X$
has $T_X\in\add T$ and factors every map from an object of $\add T$
to $X$; a left approximation is defined dually.
\end{dfn*}

\begin{dfn*}[{\cite[Section 2]{AO}; \cite[Definition 3.1]{AOTree}}]\label{def:cluster-type}
Two finite-dimensional $\ta_2$-finite algebras of global dimension
at most two are cluster equivalent if their generalized cluster
categories are $k$-linearly triangle equivalent. Such an algebra
$\La$ has acyclic cluster type $Q$ if $Q$ is a finite acyclic quiver
and $\calC_\La\simeq\calC_{kQ}$.
\end{dfn*}

Cluster equivalence need not preserve the canonical objects. Tree
cluster type means that the underlying graph of $Q$ is a tree.
An algebra is piecewise hereditary if its bounded derived category
is equivalent to that of a hereditary abelian category; see
\cite[Section 6]{AOTree}. For further background see \cite{AmiotSurvey}.

\subsection{Quivers with potential and silting objects}\label{subsec:potentials}
Let $Q$ be a finite quiver, $l=\prod_{i\in Q_0}ke_i$, and
$\mathfrak m$ its arrow ideal. Paths are multiplied in travel order.
Write $\widehat{kQ}=\varprojlim_n kQ/\mathfrak m^n$.

\begin{dfn*}[{\cite[Definitions 3.1--3.2 and 4.2]{DWZ}}]\label{def:potential}
A potential $W$ is a possibly infinite linear combination of cycles,
considered modulo the closure of the span of differences of cyclic
rotations. In this paper potentials have order at least two.
The cyclic derivative is the continuous extension of
\[
 \partial_a(a_1\cdots a_d)
 =\sum_{a_j=a}a_{j+1}\cdots a_da_1\cdots a_{j-1}.
\]
The completed Jacobian algebra is
$\Jac(Q,W)=\widehat{kQ}/\overline{(\partial_aW:a\in Q_1)}$.
The potential is Jacobi-finite if this algebra is finite-dimensional,
and reduced if its quadratic part vanishes. A right equivalence is a
continuous vertex-fixing isomorphism of completed path algebras taking
one potential to a cyclically equivalent potential.
\end{dfn*}

For Jacobi-finite $W$, put $J=\Jac(Q,W)$ and
$[W]\in HH_0(J)=J/[J,J]$, where $[J,J]$ is the linear span of
commutators. The recognition argument uses this class in addition to $J$.

\begin{dfn*}[{\cite[Section 2.6]{KY}; \cite[Section 6.2]{Keller}}]\label{def:ginzburg}
Add to $Q$ a reverse arrow $a^*:j\to i$ for every $a:i\to j$ and
a loop $t_i$ at each vertex. The completed Ginzburg dg algebra is
the degreewise arrow completion of the resulting graded path algebra,
with continuous differential determined by
\[
 \begin{gathered}
 |a|=0,\quad |a^*|=-1,\quad |t_i|=-2,\\
 d(a)=0,\quad d(a^*)=\partial_aW,\quad
 d(t_i)=e_i\bigl(\sum_a[a,a^*]\bigr)e_i.
 \end{gathered}
\]
Here the brackets are graded commutators. We denote this algebra by
$\widehat\Ga(Q,W)$. For a finite potential, the same formulas on
the uncompleted path algebra define $\Ga(Q,W)$.
\end{dfn*}

We follow Keller--Yang's cohomological convention and arrow-ideal
completion. For the ordinary model, Keller's definition in
\cite[Section 6.2]{Keller}, with $n=3$, agrees after replacing his
loop generator by its negative. The construction originates in
\cite{Ginzburg}; the distinction between the completions is explained
in \cite[Remark 2.7]{KY}.
The identity $\sum_a[a,\partial_aW]=0$ implies $d^2=0$, and
$H^0\widehat\Ga(Q,W)=\Jac(Q,W)$. For the ordinary model, $H^0$
is the ordinary derivative quotient; comparison with the completed
one requires a separate argument.

The relation-potential construction of \cite[Section 6.9]{Keller}
will be used repeatedly. Let $\La=kQ/I$ be finite-dimensional basic,
with $I$ admissible and $\operatorname{gl.dim}\La\le2$.
Choose vertexwise representatives $r$ of a basis of
$I/(\mathfrak m I+I\mathfrak m)$. Add an arrow $\ro_r$ from the
terminal vertex of $r$ to its initial vertex, obtaining $\overline Q$,
and set
\[
 W_\La=\sum_r\ro_r r,\qquad
 \deg a=0\ (a\in Q_1),\quad \deg\ro_r=1.
\]
By \cite[Theorem 6.10]{Keller}, we obtain
$\Pi_3(\La)\simeq\Ga(\overline Q,W_\La)$.
The internal grading extends by
$\deg a^*=1-\deg a$ and $\deg t_i=1$, so $\deg W_\La=1$
and $\deg d=0$.

\begin{dfn*}[{\cite[Definitions 2.1 and 2.10]{AI}}]\label{def:silting}
An object $M\in\perf R$ is presilting if
$\Hom(M,M[i])=0$ for $i>0$. It is silting if, in addition,
$\thick(M)=\perf R$, and tilting if it generates and the same
vanishing holds for every $i\ne0$. For silting objects put
$P\ge Q$ if $\Hom(P,Q[i])=0$ for $i>0$.
\end{dfn*}

For subcategories $\mathcal U,\mathcal V$, the notation
$\mathcal U*\mathcal V$ means the objects in triangles
$U\to X\to V\to U[1]$. For nonpositive $R$, two-term objects are
those in $\add R*\add R[1]$. For silting $M$ this is equivalent
to $R\ge M\ge R[1]$ \cite[Proposition 2.23]{AI}.
Reduction to $H^0R$ induces the two-term silting correspondence
\cite[Proposition 3.2(a)]{XYZ}; it is not asserted to be fully
faithful on all two-term objects. We use the cluster-category lifting
result \cite[Theorem 5.2]{HI}.

\subsection{Coverings and comparison lemmas}\label{subsec:silting-coverings}

\begin{dfn*}[{Compare \cite[Section 2.b]{AP}}]\label{def:covering}
For an internally $\bbZ$-graded dg algebra $L$, its covering dg
category has objects $\bbZ$ and morphism complexes
\[
 \Cov(L)(m,n)=\Hom^\bullet_{L,\mathrm{gr}}(L(m),L(n))\iso L_{n-m}.
\]
Composition is composition of maps of right modules. For a bigraded
algebra $R$, $\Cov_v(R)(m,n)=R_{v=n-m}$ retains the other internal
degree $u$.
\end{dfn*}

The correspondence $X\mapsto[m\mapsto X_{-m}]$ identifies graded
dg modules with right modules over this covering, and induces the
corresponding derived equivalence. For the covering quotients used
below, see \cite[Appendix A, Setting A.5 and Theorem A.7]{Han}.
We use the notation
\[
 \mathcal Q(S)=\perf^{\bbZ^2}S/
       (\perf^{\bbZ^2}S\cap\calD^{\bbZ^2}_{\fd}S).
\]
Under the hypotheses established in Section~\ref{sec:cover}, the final
comparison has the form
\[
 \Db(\gr^qL_1)\simeq\mathcal Q(A)\simeq\mathcal Q(R)
       \simeq\Db(\gr^pL_2),
\]
where $L_1=H^0(A_{p=0})$ has the residual $q$-grading and
$L_2=H^0(R_{q=0})$ has the residual $p$-grading.
Thus $\gr^qL_1$ and $\gr^pL_2$ use the two different residual gradings.
The middle equivalence comes from a bigraded Morita bimodule; the
outer equivalences require Definition~\ref{def:admissible}.

\label{subsec:conventions}
The following comparison lemmas will be used below. All internal gradings
are direct sums. None of the lemmas identifies a completion with an
ordinary graded algebra.

\begin{lem}\label{lem:forget}
Let $R$ be an ordinary bigraded dg algebra and $F$ forget the second grading.
\begin{enumerate}[label=(\roman*)]
\item If $FX$ is compact in $\calD^{\bbZ}_pR$, then $X$ is compact in $\calD^{\bbZ^2}R$.
\item If $X$ is perfect, then for any $Y$,
\[
 \Hom_{\calD^{\bbZ}_pR}(FX,FY(a)[i])
 \iso\Ds_{b\in\bbZ}
 \Hom_{\calD^{\bbZ^2}R}(X,Y(a,b)[i]).
\]
\item If $X$ is perfect, then, after forgetting the second grading, $\uREnd^{\bbZ^2}_R(X)$ and $\uREnd^p_R(FX)$ are quasi-isomorphic as $p$-graded dg algebras.
\end{enumerate}
The same assertions hold for grading blocks $p\in\bbZ^r$, $q\in\bbZ^s$, and when both blocks are forgotten.
\end{lem}
\begin{proof}
For a $p$-graded map $f$ and a $q$-homogeneous element $x\in X_q$, put
\[
 f_b(x)=\operatorname{pr}_{Y_{q+b}}f(x).
\]
Each $f_b$ is $R$-linear and $(d_{\Hom}f)_b=d_{\Hom}(f_b)$. In particular, a degree-zero map which is null-homotopic after forgetting $q$ already has a degree-zero null-homotopy. We may apply this observation to derived morphisms by choosing graded semi-free resolutions; these remain semi-free after forgetting $q$.

(i) Let $f:X\to\Ds_jY_j$. Since $FX$ is compact, there is a finite set $I$ such that $Ff$ factors through $\Ds_{j\in I}FY_j$. Taking the $q$-degree-zero component, we obtain
\[
 X\longrightarrow\Ds_{j\in I}Y_j
 \longrightarrow\Ds_jY_j
\]
in the graded derived category. The same operation on a null-homotopy shows that no additional relations appear after forgetting $q$. Thus $\Hom(X,-)$ commutes with coproducts.

(ii) For a shift of $R$, the formula follows from the homogeneous decomposition of $H^*(Y)$. Both sides preserve finite sums, triangles and direct summands in $X$. The formula therefore holds for every perfect $X$.

(iii) Choose a graded semi-free resolution $Q$ of $X$. Finite sums of homogeneous endomorphisms form a dg subalgebra of the endomorphism algebra after forgetting $q$:
\[
 F\underline{\End}^{\bbZ^2}_R(Q)\longrightarrow\underline{\End}^p_R(FQ).
\]
By (ii), this map induces an isomorphism on cohomology in each $p$-degree. Hence it is a quasi-isomorphism. The same proof applies to grading blocks, or to forgetting all internal degrees.
\end{proof}

\begin{lem}\label{lem:idempotents}
Let $k$ be algebraically closed. A finite-dimensional basic $\bbZ^r$-graded ordinary $k$-algebra $J$ has a complete set of primitive orthogonal idempotents of internal degree zero.
\end{lem}
\begin{proof}
Let $T=(\mathbb G_{m,k})^r$ act on $J$ through the grading. This action is rational and preserves $\rad J$. Connectedness implies that $T$ acts trivially on $J/\rad J\iso k^n$, since it cannot permute its finitely many primitive idempotents nontrivially. Since $T$ is linearly reductive, we obtain the following surjection by taking $T$-invariants:
\[
 0\longrightarrow(\rad J)^T\longrightarrow J^T
 \longrightarrow k^n\longrightarrow0,
 \qquad ((\rad J)^T)^N=0\quad(N\gg0).
\]
The primitive idempotents of the quotient lift to elements $e_i\in J^T=J_0$ satisfying
\[
 e_ie_j=\de_{ij}e_i,\qquad \sum_i e_i=1,
 \qquad \overline e_i=(0,\ldots,1,\ldots,0).
\] They are primitive in $J$ since their images are primitive in its semisimple quotient.
\end{proof}

\begin{lem}\label{lem:proj}
Let $R$ be a nonpositive dg algebra with internal grading $p\in\bbZ^r$, put $J=H^0R$, and
$\mathcal P_p=\add\{R(a):a\in\bbZ^r\}\subset\calD^{\bbZ^r}_pR$.
Then $H^0$ identifies $\mathcal P_p$ with the finitely generated graded projective $J$-modules. For $P\in\mathcal P_p$ and $X\in\calD^{\bbZ^r}_pR$,
\[
 \Hom_{\calD^{\bbZ^r}_pR}(P,X)
 \iso\Hom_{\operatorname{Gr}_pJ}(H^0P,H^0X).
\]
Here $\operatorname{Gr}_pJ$ denotes all graded modules; $H^0X$ need not be finitely generated. The same assertions hold for any finite-rank internal grading group.
\end{lem}
\begin{proof}
For $P=R(a)$ the formula follows from the definition of cohomology. Both sides preserve finite sums and direct summands, so the formula holds for every $P\in\mathcal P_p$. Conversely, every finitely generated graded projective is the image of an idempotent on a finite graded free module. Full faithfulness lifts this idempotent, and it splits in the derived category.
\end{proof}

\section{Relative grading of two-term silting objects}\label{sec:relative}
In this section we prove that a two-term silting object admits an additional internal grading without changing its given graded realization. We begin with bounded complexes of projective modules over a finite-dimensional algebra.

Huisgen-Zimmermann and Saor\'in relate infinitesimal deformations of projective complexes to self-extensions and prove openness of rigid orbits \cite[Theorem~7 and Corollary~8]{HZS}. Here the obstruction lies in the homotopy category preserving the first grading. Its vanishing permits a correction of the Euler operator.

\subsection{Gradings on complexes of projective modules}
The following proposition allows us to extend a grading on a rigid complex.
\begin{prp}\label{prop:finite}
Let $k$ be algebraically closed of characteristic zero, and let $J$ be a finite-dimensional $(\bbZ^r\ds\bbZ^s)$-graded $k$-algebra, with a complete set of primitive orthogonal idempotents in degree zero. Let
\[
 P\in \Kb(\gr^p\prj J),\qquad
 \Hom_{\Kb(\gr^p\prj J)}(P,P[1])=0,
 \qquad p\in\bbZ^r.
\]
Then each term of $P$ admits a $(\bbZ^r\ds\bbZ^s)$-grading extending its given $p$-grading such that the terms are graded projective and the differential is homogeneous of degree zero.
\end{prp}
\begin{proof}
Suppose first that $s=1$.

Each $P^i$ is a finite sum of $p$-graded shifts of $eJ$, where $e$ is one of the given primitive idempotents. To see this, lift homogeneous generators of the semisimple top; the decomposition follows from projectivity and Nakayama's lemma. Here $\rad J$ is homogeneous because the grading torus preserves it.

Assign $q$-degree zero to the generator of each summand $eJ$. Define $k$-linear maps $E_J$ and $E$ by
\[
 E_J(a)=ma\quad(a\in J_{*,m}),\qquad
 E(x)=nx\quad(x\text{ of temporary }q\text{-degree }n).
\]
Thus $E$ preserves $p$-degree and satisfies $E(xa)=E(x)a+xE_J(a)$. It need not commute with $d$.

Put
\[
 \mathcal E^n=\prod_i\Hom_{\mathrm{gr}_pJ}(P^i,P^{i+n}),\qquad
 \partial(f)=df-(-1)^nfd.
\]
Since $P$ is bounded, these products are finite. Moreover,
\[
 H^1(\mathcal E)=\Hom_{\Kb(\gr^p\prj J)}(P,P[1])=0.
\]
Put $c=Ed-dE$. By the Leibniz identity, $c\in\mathcal E^1$, and
\[
 \partial(c)=dc+cd=0,\qquad
 H^1(\mathcal E)=0
 \ \Longrightarrow\ \exists h\in\mathcal E^0:
 c=\partial(h)=dh-hd.
\]
Put $\nabla=E+h$. Then
\begin{equation}\label{eq:connection}
 d\nabla=\nabla d,\qquad
 \nabla(xa)=\nabla(x)a+xE_J(a).
\end{equation}

Since $k$ is algebraically closed and $\Ds_iP^i$ is finite dimensional, it decomposes into generalized eigenspaces of $\nabla$. Write $P^i(\la)$ for the generalized $\la$-eigenspace. By \eqref{eq:connection}, we have
\[
 dP^i(\la)\subseteq P^{i+1}(\la),\qquad
 P^i(\la)J_{*,m}\subseteq P^i(\la+m).
\]
For the second inclusion, iterate
$(\nabla-\la-m)(xa)=((\nabla-\la)x)a$.

Choose a set $L$ of representatives for the classes of eigenvalues
in $k/\bbZ$, and put
\[
 P^i_{p,n}=\bigoplus_{\la_0\in L}
       \bigl(P^i_p\cap P^i(\la_0+n)\bigr).
\]
Since $\operatorname{char}k=0$, each eigenvalue occurs for a unique pair
$(\la_0,n)$. The operator $\nabla$ preserves $p$-degree. Therefore
\[
 P^i_p=\Ds_nP^i_{p,n},\qquad
 P^i_{p,n}J_{u,m}\subseteq P^i_{p+u,n+m},\qquad
 dP^i_{p,n}\subseteq P^{i+1}_{p,n}.
\]
It remains to check graded projectivity. A graded free surjection onto $P^i$ has an ungraded section. Its multidegree-zero component is a graded section, as required.

For $s>1$, we argue by induction. By projecting a null-homotopy to degree zero in the added coordinates, we obtain
\[
 \Hom_{\Kb(\gr_{\bbZ^{r+j}}\prj J)}(P_j,P_j[1])
 \incl
 \Hom_{\Kb(\gr^p\prj J)}(P,P[1])=0.
\]
The case $s=1$ therefore applies at each stage, without changing the preceding gradings or the differential.
\end{proof}

\begin{rmk}
Write $F$ for forgetting $p$. The map
\[
 \Hom_{\Kb(\gr^p\prj J)}(P,P[1])
 \longrightarrow\Hom_{\Kb(\prj J)}(FP,FP[1])
\]
is injective, by projection of null-homotopies to $p$-degree zero. The hypothesis is therefore weaker than ungraded rigidity. For one additional coordinate, the proof only requires the particular map $Ed-dE:P\to P[1]$ to be null-homotopic in the $p$-graded category.
\end{rmk}

\subsection{Gradings on two-term silting objects}
We next apply \prpref{prop:finite} to dg algebras. We use the two-term silting correspondence recalled in Subsection~\ref{subsec:potentials}.
\begin{lem}\label{lem:presentation}
Let $k$ be algebraically closed of characteristic zero. Let $R$ be nonpositive and $p$-graded, with $p\in\bbZ^r$ and $J=H^0R$ finite dimensional and basic. If $N\in\perf^{\bbZ^r}_pR$ is two-term after forgetting $p$, then there is a triangle
\[
 P_1\ya{d}P_0\longrightarrow N\longrightarrow P_1[1]
\]
in $\calD^{\bbZ^r}_pR$ with $P_0,P_1\in\mathcal P_p$.
\end{lem}
\begin{proof}
The finite-dimensional module $H^0N$ has a graded projective cover
$\overline P_0\to H^0N$. Here graded covers are ordinary covers after forgetting: the radical of $J$ is homogeneous and the primitive idempotents can be chosen in degree zero. Lift the cover to $f:P_0\to N$ by \lemref{lem:proj}, and define $P_1$ by its cone triangle.

After forgetting $p$, choose any two-term presentation
$Q_1\to Q_0\ya{g}N\to Q_1[1]$ with $Q_i\in\add R$.
Both $f$ and $g$ are right $\add R$-approximations. For $f$ this follows from its surjectivity on $H^0$; for $g$ it follows from $\Hom(R,Q_1[1])=0$. The map $f$ is right minimal: an endomorphism $u$ with $fu=f$ induces an automorphism on $H^0P_0$ by the cover property, and hence is an automorphism by full faithfulness on $\add R$.

Choose $a:P_0\to Q_0$ and $b:Q_0\to P_0$ such that
\[
 f=ga,
 \qquad g=fb,
 \qquad f=fba.
\]
By right minimality, $ba\in\operatorname{Aut}(P_0)$. Hence $a$ is a split monomorphism. With respect to the induced decomposition, we have
\[
 Q_0\iso P_0\ds Q',\quad g=(f,0),\quad
 Q_1\iso P_1\ds Q',\quad P_1\in\add R.
\]

Consequently $H^0P_1$ is projective after forgetting $p$. It is graded projective: a graded free surjection splits ungraded, and the degree-zero part of the section is a graded section. Lift $H^0P_1$ to $P'_1\in\mathcal P_p$ and its identity map to $P'_1\to P_1$ by \lemref{lem:proj}. After forgetting $p$, both objects belong to $\add R$, and this map is an isomorphism since its $H^0$ is an isomorphism. Forgetting reflects quasi-isomorphisms, so it is also an isomorphism in $\calD^{\bbZ^r}_pR$. Replace $P_1$ by $P'_1$.
\end{proof}

\begin{thm}\label{thm:relative}
Let $k$ be algebraically closed of characteristic zero, and let $R$ be a nonpositive $(\bbZ^r\ds\bbZ^s)$-graded dg $k$-algebra with $H^0R$ finite dimensional and basic. Write $p\in\bbZ^r$ and $q\in\bbZ^s$ for the two grading blocks. Let $N\in\perf^{\bbZ^r}_pR$ have underlying ungraded object two-term silting. Then there exists $\widetilde N\in\perf^{\bbZ^{r+s}}R$ such that
\begin{equation}\label{eq:relative-target}
 \Fgt_q\widetilde N\simeq N\quad\text{in }\perf^{\bbZ^r}_pR.
\end{equation}
Moreover
\[
 \Fgt_q\uREnd^{\bbZ^{r+s}}_R(\widetilde N)\simeq\uREnd^p_R(N)
\]
as $p$-graded dg algebras. The internal shifts of $\widetilde N$ form a family of compact generators of $\calD^{\bbZ^{r+s}}R$.
\end{thm}
\begin{proof}
Put $J=H^0R$. By \lemref{lem:presentation}, there is a triangle in $\calD^{\bbZ^r}_pR$
\[
 P_1\ya{d}P_0\longrightarrow N\longrightarrow P_1[1].
\]
Its image under $-\Lox_RJ$ is the two-term complex
$P=(H^0P_1\to H^0P_0)$. After forgetting $p$, this is two-term silting by \cite[Proposition 3.2(a)]{XYZ}. In particular $\Hom(P,P[1])=0$. By \lemref{lem:idempotents}, $J$ has a complete set of primitive orthogonal idempotents of bidegree zero. Thus \prpref{prop:finite} applies to $P$ and extends its $p$-grading to a $(p,q)$-grading without changing the differential.

Use the multigraded version of \lemref{lem:proj} to lift its projective terms to $\widetilde P_1,\widetilde P_0$ and its differential to a degree-zero morphism $\widetilde d$. Define $\widetilde N$ as its cone. After forgetting $q$, the terms are identified with $P_1,P_0$ by their graded $H^0$. Under these identifications the lifted differential equals $d$ in the $p$-graded derived category, again by full faithfulness on the projective terms. By comparing the two cone triangles, we obtain \eqref{eq:relative-target}.

By \lemref{lem:forget} and \eqref{eq:relative-target},
\[
 \Fgt_q\uREnd_R^{\bbZ^{r+s}}(\widetilde N)
 \simeq\uREnd_R^p(\Fgt_q\widetilde N)
 \simeq\uREnd_R^p(N).
\]
To prove generation, we apply \lemref{lem:forget} and obtain
\[
 \begin{aligned}
 &\Hom_{\calD^{\bbZ^{r+s}}R}(\widetilde N(a,b),X[i])=0
       &&(a\in\bbZ^r,\ b\in\bbZ^s,\ i\in\bbZ)\\
 &\hspace{12mm}\Longrightarrow\ \Hom_{\calD R}(N,\Fgt X[i])=0
       &&(i\in\bbZ)\\
 &\hspace{12mm}\Longrightarrow\ \Fgt X=0\ \Longrightarrow\ X=0.
 \end{aligned}
\]
Since these shifts are compact, it follows from \cite[Theorem 3.4]{KellerDG} that
\[
 \thick\{\widetilde N(a,b)\}=\perf^{\bbZ^{r+s}}R.
\]
\end{proof}

\begin{cor}\label{cor:replacement}
Suppose that $R'$ and $R$ are $p$-graded dg algebras connected by a zigzag of graded quasi-isomorphisms. Let $N$ be the image of $N'\in\perf^{\bbZ^r}_pR'$ under the induced equivalence. If the underlying $N'$ is two-term silting relative to $R'$, the object $N$ is two-term silting relative to $R$. Whenever $R$ satisfies \thmref{thm:relative}, its bigraded lift has $p$-graded endomorphism dg algebra quasi-isomorphic to that of the original $N'$.
\end{cor}
\begin{proof}
Extension and restriction along dg quasi-isomorphisms identify regular objects. They therefore preserve the two-term interval, silting and the specified $p$-graded object. By \thmref{thm:relative}, the required lift exists. The last assertion follows by comparing the endomorphism dg algebras under these quasi-equivalences.
\end{proof}

\begin{rmk}
Taking $r=0$ and $s=1$ in \thmref{thm:relative}, every two-term silting object admits a grading for the given internal grading on $R$. The two-term hypothesis is used when lifting the single differential from $H^0R$: for longer dg objects, the composites and their higher homotopies must also be lifted.
\end{rmk}

\section{A common multigraded dg model}\label{sec:common}
Throughout this section, $k$ is an uncountable algebraically closed field of characteristic zero, and $r,s\ge0$. We construct a common multigraded model for two internally graded dg algebras which are quasi-isomorphic after forgetting their gradings.

The two gradings induce torus actions on the bar coalgebra of a common minimal $\Ainf$ model. We first conjugate one action to make the actions commute, and then linearize them. The resulting action preserves the full coderivation, not only the multiplication on cohomology. This construction is separate from the ordinary potential recognition in Section~\ref{sec:completion}.

We shall use the following elementary lemma to choose compatible solutions on the finite-dimensional truncations.
\subsection{Compatible actions on the bar coalgebra}
\begin{lem}\label{lem:countable}
Let $k$ be an uncountable algebraically closed field, and let $K\subset k$ be a countable subfield. Let $\mathcal E$ be a countable family of polynomial equations over $K$ in countably many variables. If every finite subfamily has a solution in $k$, then the entire family has a solution in $k$.
\end{lem}
\begin{proof}
Let $I\subset K[x_1,x_2,\ldots]$ be the ideal generated by the equations. Finite solvability implies
\[
 1=\sum_{j=1}^m a_j f_j\quad(f_j\in\mathcal E)
 \ \Longrightarrow\ 1=0
\]
after evaluation at a common zero of $f_1,\ldots,f_m$. Hence $I\ne(1)$. Choose a maximal ideal $\mathfrak n\supset I$ and set
\[
 L=K[x_1,x_2,\ldots]/\mathfrak n,\qquad |L|\le\aleph_0.
\] The extension $L/K$ embeds in $k$ over $K$: choose images of a transcendence basis algebraically independent over $K$ and extend the embedding across the remaining algebraic extension. There is enough transcendence degree because $K$ is countable and $k$ is uncountable. The images of the residue classes of the $x_i$ solve every equation.
\end{proof}

\begin{prp}\label{prop:torus}
Let $k$ be an uncountable algebraically closed field of characteristic zero, and let $r,s\geq0$. Let $V=\Ds_{j<0}V^j$ be a negatively cohomologically graded $k$-vector space with each $V^j$ finite dimensional. Let
\[
 C=T^c(V)=\Ds_{n\geq0}V^{\ox n}
\]
be the coaugmented tensor coalgebra and let $b$ be a degree-one square-zero coderivation with zero constant term. Suppose that $\ro$ and $\si$ are degree-preserving coalgebra actions of split tori $T_p=(\mathbb G_{m,k})^r$ and $T_q=(\mathbb G_{m,k})^s$ on $C$, commuting with $b$, rational on every cohomological degree. Then there is a degree-preserving coalgebra automorphism $g$, commuting with $b$, such that $\ro$ and $g\si g^{-1}$ commute.

The same assertion holds for strictly unital $\Ainf$ automorphisms, with the unit fixed, when $b$ encodes a strictly unital structure.
\end{prp}

\begin{proof}
We divide the proof into four steps.

\medskip\noindent\emph{Step 1.}
For $N\ge0$, put
\[
 C_N=\Ds_{-N\le j\le0}C^j
 =k\ds\Ds_{n=1}^{N}
 \ \Ds_{\substack{j_1,\ldots,j_n<0\\-N\le j_1+\cdots+j_n<0}}
 V^{j_1}\ox\cdots\ox V^{j_n}.
\]
Since each $j_i<0$, the inequalities $-N\le j_1+\cdots+j_n$ imply $n\le N$ and $-N\le j_i<0$. Hence
\[
 \dim_k C_N<\infty,\qquad
 \De(C_N)\subset C_N\ox C_N,\qquad
 b(C_N)\subset C_N,\qquad C=\bigcup_{N\ge0}C_N.
\]
The inclusion involving $b$ follows from $|b|=1$ and $b(1)=0$. Choose compatible homogeneous bases of these spaces. Write $\De_N,b_N,\epsilon_N$ for the restricted coproduct, differential and counit, and $\io_N:C_N\incl C_{N+1}$ for the inclusion.

\medskip\noindent\emph{Step 2.}
Let $G_N$ consist of the degree-preserving invertible matrices $X$ satisfying
\[
 \De_NX=(X\ox X)\De_N,\quad
 b_NX=Xb_N,\quad X(1)=1,\quad\epsilon_NX=\epsilon_N.
\]
These are polynomial identities. Together with a matrix variable $Y$, the equations
\[
 XY=YX=I_{C_N}
\]
realize $G_N$ as an affine algebraic group. In the strictly unital case we also impose
\[
 X^{(1)}(e)=e,\qquad
 X^{(n)}(v_1,\ldots,v_n)=0
 \quad(n\ge2,\ v_i=e\text{ for some }i),
\]
whenever the input belongs to $C_N$, where $e$ is the suspended unit and $X^{(n)}$ is the Taylor component obtained by projection to $V$. These are polynomial equations; composition and inversion preserve them. We use this subgroup in the unital case.

Let $\ro_N,\si_N$ be the restricted actions. Their images are tori in $G_N^\circ$. Choose maximal tori $T_N,U_N\subset G_N^\circ$ containing these images, respectively. By conjugacy of maximal tori, we have
\[
 z_NU_Nz_N^{-1}=T_N
 \quad\text{for some }z_N\in G_N^\circ;
\]
see \cite[Expos\'e XII, Theorem 6.6(a)]{SGA}. Thus
\[
 \ro_N(T_p),\ z_N\si_N(T_q)z_N^{-1}\subset T_N,
 \qquad
 [\ro_N(t),z_N\si_N(u)z_N^{-1}]=0.
\]
Here and below brackets in a matrix commutation equation denote the additive commutator. The elements $z_N$ are not yet chosen compatibly.

\medskip\noindent\emph{Step 3.}
Introduce degree-preserving matrix variables $X_N,Y_N$ for every $N$. Require $X_N\in G_N$ and impose
\[
 \begin{aligned}
 X_NY_N&=Y_NX_N=I_{C_N},\\
 X_{N+1}\io_N&=\io_NX_N,\qquad
 Y_{N+1}\io_N=\io_NY_N.
 \end{aligned}
\]
Since the actions are rational, we may write
\[
 \ro_N(t)=\sum_{\al\in A_N}t^\al R_{N,\al},\qquad
 \si_N(u)=\sum_{\be\in B_N}u^\be S_{N,\be},
\]
where $A_N\subset\bbZ^r$ and $B_N\subset\bbZ^s$ are finite. The desired commutation identity is therefore equivalent to the finite family
\[
 R_{N,\al}X_NS_{N,\be}Y_N
 =X_NS_{N,\be}Y_NR_{N,\al}
 \qquad(\al\in A_N,\ \be\in B_N).
\]
These identities are obtained by comparing the coefficients of $t^\al u^\be$. Together with the equations defining $G_N$, they form a countable polynomial system $\mathcal E$. Its coefficients generate a countable subfield $K\subset k$.

\medskip\noindent\emph{Step 4.}
Let $\mathcal E_0\subset\mathcal E$ be finite, and choose $M$ containing every index occurring in $\mathcal E_0$. For $N\le M$, restriction of degree-preserving automorphisms defines homomorphisms
\[
 r_{M,N}:G_M\longrightarrow G_N,\qquad
 r_{L,N}r_{M,L}=r_{M,N}\quad(N\le L\le M).
\]
By Step 2, there exists $z_M\in G_M^\circ$ with the stated conjugacy property. Put
\[
 X_N=r_{M,N}(z_M),\qquad
 Y_N=r_{M,N}(z_M^{-1})\qquad(N\le M).
\]
These matrices solve $\mathcal E_0$: the group identities and commutation identities restrict from $C_M$, and the restriction equations hold by construction. No surjectivity of $r_{M,N}$ is required.

By \lemref{lem:countable}, the system $\mathcal E$ has a solution $(X_N,Y_N)_{N\ge0}$. Define
\[
 g(x)=X_Nx,\qquad h(x)=Y_Nx\qquad(x\in C_N).
\]
The restriction equations make these definitions independent of $N$, and the remaining equations give
\[
 gh=hg=1_C,\qquad
 \De g=(g\ox g)\De,\qquad bg=gb,
\]
\[
 \ro(t)g\si(u)h=g\si(u)h\ro(t).
\]
Thus $h=g^{-1}$ and $g$ has the required properties. In the unital case, every Taylor component satisfies the imposed unit equations at a sufficiently large finite level, so $g$ is strictly unital as well.
\end{proof}

\begin{rmk}\label{rem:field}
The argument requires finite solvability of the full system, not merely nonemptiness at each finite level. Uncountability is used in \lemref{lem:countable}, where the countable residue field is embedded in $k$ over the coefficient field.
\end{rmk}

\subsection{Commuting actions and minimal models}
\begin{lem}\label{lem:linear}
In the setting of \prpref{prop:torus}, a rational torus action on $C$ is conjugate, by a coalgebra automorphism, to the tensor extension of its action on the primitive subspace $V$. If one torus factor already acts by a tensor extension, the conjugating automorphism can be chosen to commute with that factor. Strict units can be preserved.
\end{lem}
\begin{proof}
The primitive subspace $V$ is stable under every coaugmented coalgebra automorphism. Degree by degree, choose a torus-equivariant retraction $q:C\to V$ of its inclusion, with $q(1)=0$. This is possible because the homogeneous pieces of $C$ are finite-dimensional rational representations and a torus is linearly reductive. By the universal property of $T^c(V)$, there is a unique coaugmented coalgebra map $\Th:C\to T^c(V)$ whose composite with the projection to $V$ is $q$. Its linear component is the identity. Put $F_nC=\Ds_{j\le n}V^{\ox j}$ and $u=\Th-1$. By the formula for morphisms of tensor coalgebras, we obtain
\[
 u(F_nC)\subset F_{n-1}C,\qquad u^{n+1}|_{F_nC}=0,\qquad
 \Th^{-1}|_{F_nC}=\sum_{j=0}^{n}(-u)^j.
\]
These finite inverses agree under $F_nC\subset F_{n+1}C$, and hence define a coalgebra inverse on $C$. It intertwines the given torus action and the tensor extension of the action on $V$.

If a factor $\ro$ already acts linearly, equivariance of $q$ with respect to that factor says $\Th\ro=\ro\Th$. In the unital situation, let $D\subset C$ be spanned by tensors of length at least two that contain the suspended unit. Strict unitality makes $D$ stable under the torus: a block containing the unit cannot be contracted by a higher component, and a tensor containing the unit cannot be contracted to length one. Moreover $D\cap V=0$. Choose the equivariant retraction on $C/D$ and pull it back. Then $q(D)=0$, so $\Th$ is strictly unital.
\end{proof}

\begin{cor}\label{cor:common}
Let $H$ be a strictly unital minimal $\Ainf$ algebra over $k$, with $H^i=0$ for $i>0$ and $\dim H^i<\infty$ for all $i$. Suppose two internal gradings are represented by rational strictly unital $\Ainf$ actions of $T_p,T_q$ on its bar coalgebra, and suppose the first action is linear. There is an $\Ainf$-isomorphic minimal model with a $(\bbZ^r\ds\bbZ^s)$-grading on which every operation $m_n$ has internal bidegree $(0,0)$. The first grading is preserved up to a first-grading-preserving $\Ainf$ isomorphism, and the second is preserved up to an $\Ainf$ isomorphism.
\end{cor}
\begin{proof}
Use the convention $H[1]^i=H^{i+1}$, so $V=H[1]$ is strictly negative. Apply \prpref{prop:torus}, then \lemref{lem:linear}. Transport the coderivation by $\Th$. It commutes with the resulting linear $T_p\times T_q$ action; equivalently, every $m_n$ preserves its two weights. The first action was unchanged by the initial conjugation of the second and is preserved by $\Th$.
\end{proof}

\begin{thm}\label{thm:common-dg}
Let $k$ be an uncountable algebraically closed field of characteristic zero, and let $r,s\geq0$. Let $(A_1,p)$ and $(A_2,q)$ be unital dg $k$-algebras graded by $p\in\bbZ^r$ and $q\in\bbZ^s$, respectively. Assume that both are connective and cohomologically locally finite. If their underlying ordinary dg algebras are quasi-isomorphic, there exists an ordinary $(\bbZ^r\ds\bbZ^s)$-graded dg algebra $R$ together with zigzags of quasi-isomorphisms of dg algebras
\[
 (R,p)\simeq(A_1,p),\qquad (R,q)\simeq(A_2,q).
\]
The algebra $R$ can be chosen nonpositive.
\end{thm}
\begin{proof}
Suppose first that $H^*(A_i)=0$. Choose $h_i\in A_i^{-1}$ of internal degree zero with $dh_i=1$; take the degree-zero component of any primitive of $1$. The acyclic dg algebra $k\langle h\rangle$, with $|h|=-1$, $dh=1$, and both internal gradings trivial, maps quasi-isomorphically to each $A_i$ by $h\mapsto h_i$. If the algebras are zero the same maps are interpreted in the category allowing the zero unital algebra. This proves the assertion in the acyclic case.

Otherwise the cohomology unit is nonzero. Choose strictly unital minimal models $H_i\to A_i$ by homogeneous splittings containing the unit for the respective internal gradings; see \cite[Sections 3.3 and 8.1]{KellerAinf} and \cite[Section 1.1]{COS}. All transferred operations and morphism components preserve that grading. The zigzag of quasi-isomorphisms between the underlying dg algebras induces a strictly unital $\Ainf$ isomorphism $F:H_2\to H_1$ between the minimal models.

Put $H=H_1$. Let $\ro$ be the first action on $T^c(H[1])$, and define the second action by $\si=F\ro_qF^{-1}$. Every cohomological component of this coalgebra is finite dimensional. The transported action is rational on each such component, since there it is conjugate by a fixed finite matrix to the rational original action.

Apply \prpref{prop:torus} and \lemref{lem:linear}. In their notation put $\si'=g\si g^{-1}$ and $b'=\Th b\Th^{-1}$. The resulting minimal algebra $H'$ is bigraded. If $\io_i:H_i\to A_i$ are the original minimal-model maps, the comparison maps are
\[
 f_1=\io_1\Th^{-1},\qquad
 f_2=\io_2F^{-1}g^{-1}\Th^{-1}.
\]
Let $\ro'_q$ denote the linear second action on $H'$. The identities
\[
 \Th\ro=\ro\Th,\qquad
 \Th\si'=\ro'_q\Th,\qquad
 \si'=gF\ro_qF^{-1}g^{-1}
\]
imply that
\[
 f_1\ro(t)=\ro_{A_1,p}(t)f_1,\qquad
 f_2\ro'_q(t)=\ro_{A_2,q}(t)f_2.
\]
They are identities of coalgebra morphisms, hence include every Taylor component, not only the induced maps on cohomology. Equivalently, for every $n\ge1$, their Taylor components satisfy
\[
 \begin{aligned}
 (f_1)_n(H'_{p_1,*},\ldots,H'_{p_n,*})
   &\subseteq(A_1)_{p_1+\cdots+p_n},\\
 (f_2)_n(H'_{*,q_1},\ldots,H'_{*,q_n})
   &\subseteq(A_2)_{q_1+\cdots+q_n}.
 \end{aligned}
\]
Their linear components induce isomorphisms on cohomology, since those of $F$, $g$, $\Th$ and the minimal-model maps do. Thus both maps are strictly unital $\Ainf$ quasi-isomorphisms with the required separate gradings.

Use the unital enveloping dg algebra $U(H')$ of \cite[Proposition 2.1]{COS}. Unlike the reduced bar-cobar construction for augmented algebras, this construction applies to general strictly unital algebras. It is made from ordinary tensor direct sums and the homogeneous unit relations. Therefore it inherits the bigrading, with differential of internal bidegree zero. Its universal property turns $f_1,f_2$ into dg algebra maps $U(H')\to A_i$ preserving $p,q$, respectively. Write $u:H'\to U(H')$ for the unit and $\widehat f_i:U(H')\to A_i$ for the resulting dg maps. Then $\widehat f_i u=f_i$. On cohomology,
\[
 H^*(\widehat f_i)\,H^*(u_1)=H^*((f_i)_1),
 \qquad H^*(u_1),\ H^*((f_i)_1)\text{ invertible};
\]
thus $\widehat f_i$ is a quasi-isomorphism. The tensor generators and unit relations are homogeneous, so these maps retain their respective internal gradings. Finally put
\[
 R=\ta_{\le0}U(H'),\qquad R\incl U(H')\text{ a bigraded quasi-isomorphism}.
\]
\end{proof}

\begin{rmk}
The two comparisons preserve the given graded quasi-isomorphism types separately, without prescribing the underlying ungraded comparison. The quasi-isomorphism hypothesis cannot be replaced here by an arbitrary dg Morita equivalence. In the main proof, the theorem is applied to the endomorphism dg algebra of a selected silting object and the second completion, not to the two original completions directly.
\end{rmk}

\subsection{Local finite parts of completed torus representations}\label{subsec:localfinite}
For later comparison with completed models, we record a property of the locally finite part of a continuous torus action. Here the ground field is arbitrary.

\begin{prp}\label{prop:localfinite}
Let $k$ be any field and let $C$ be a degreewise pseudocompact dg $k$-algebra with a continuous action of a split torus $T=(\mathbb G_{m,k})^r$. Assume there is a basis of $T$-stable open dg ideals whose finite-dimensional quotients are rational $T$-representations, and that every $H^iC$ is finite dimensional. Put
\[
 C_{\mathrm{lf}}=\Ds_{\chi\in X^*(T)}C_\chi.
\]
Then $C_{\mathrm{lf}}\incl C$ is an ordinary dg quasi-isomorphism.
\end{prp}
\begin{proof}
Weight projections on the finite-dimensional quotients commute with the transition maps and define continuous projections $e_\chi:C\to C_\chi$. They commute with $d$. Multiplication adds weights, so $C_{\mathrm{lf}}$ is an ordinary graded dg subalgebra.

The image of a continuous map of linearly compact vector spaces is closed. Thus the cohomology of $C$ is Hausdorff pseudocompact; this is also the mechanism of \cite[Lemma 7.12]{KY}. Finite-dimensionality makes each $H^iC$ discrete. The induced torus representation is rational because its quotient map factors through a finite-dimensional rational quotient of the cycle space. It therefore has only finitely many weights. For a cohomology class $\xi=\sum_{\chi\in F}\xi_\chi$ with $F$ finite and a representative $z$, put
\[
 z_F=\sum_{\chi\in F}e_\chi z\in C_{\mathrm{lf}},\qquad
 dz_F=0,\qquad[z_F]=\sum_{\chi\in F}e_\chi\xi=\xi.
\]
This proves surjectivity on cohomology.

If a finite-weight cycle $z=\sum_{\chi\in F}z_\chi$ is $dy$ in $C$, then
\[
 z=d\left(\sum_{\chi\in F}e_\chi y\right)
\]
already in $C_{\mathrm{lf}}$. This proves injectivity.
\end{proof}

\begin{cor}\label{cor:localfinite}
If $C$ has commuting actions of $T_p,T_q$ satisfying the preceding hypotheses jointly, its joint local finite part $C^{p,q}$ maps quasi-isomorphically to each of its single local finite parts $C^p,C^q$. The first map preserves $p$ and the second preserves $q$.
\end{cor}
\begin{proof}
By \prpref{prop:localfinite}, we obtain
\[
 H^*(C^{p,q})\ya{\sim}H^*(C),\qquad
 H^*(C^p)\ya{\sim}H^*(C),\qquad
 H^*(C^q)\ya{\sim}H^*(C).
\]
The factorizations $C^{p,q}\incl C^p\incl C$ and $C^{p,q}\incl C^q\incl C$ prove the claim by two-out-of-three.
\end{proof}

\begin{rmk}\label{rem:formal-weights}
A pure weight space can contain infinite formal sums of paths. The local finite condition bounds the number of weights of an element, not the number of its paths. Consequently this result does not identify the local finite part with the polynomial path algebra.
\end{rmk}

\section{Ordinary and completed dg models}\label{sec:completion}
In this section we compare the ordinary Calabi--Yau completion with the arrow-completed Ginzburg algebra. For a specified silting object, we construct a continuous Morita bimodule and apply potential recognition to its endomorphism dg algebra. We thereby obtain the ordinary quasi-isomorphism required in \thmref{thm:common-dg}.

All completions in this section are taken degree by degree in the cohomological grading. They are not completions of the total ungraded vector space.

\subsection{Ordinary and arrow-completed tensor models}
\begin{lem}\label{lem:completion}
Let $l=\bbC^r$ and $A=(T_lV,d)$, with $V$ finite dimensional and concentrated in nonpositive degrees, and $d(V)\subset(V)^2$. Let
$B=\widehat A=\varprojlim_nA/(V)^n$.
Assume that $A$ and $B$ are cohomologically locally finite and that their augmentations identify $l$ with the full semisimple quotients of $H^0A$ and $H^0B$. Then the canonical map $A\to B$ is an ordinary dg quasi-isomorphism.
\end{lem}
\begin{proof}
Let $I_A,I_B$ be the augmentation ideals. By decomposing each nontrivial path into its first arrow and the remaining path, we obtain
\[
 I_A=V\ox_lA,\qquad I_B=V\ox_lB
\]
as graded right modules. The second equality uses finite-dimensionality of $V$. Equip the right-hand sides with the induced differentials. If $v$ has degree $j$, the first arrow of every term of $d(v)$ has degree at least $j+1$, since all remaining arrow degrees are nonpositive. Ordering generators by decreasing degree shows that these modules are finite semi-free.

Thus $P_A=\Cone(I_A\to A)$ and $P_B=\Cone(I_B\to B)$ are finite semi-free resolutions of $l$, with
\[
 P_A\ox_AB\iso P_B.
\]
For $F=-\Lox_AB$, we obtain $F(l)\simeq l$. Derived adjunction, with the augmentation identification on restriction, implies that the map on all self-extensions of $l$ induced by $F$ is an isomorphism. Consequently
\[
 E_A=\End_A(P_A)\longrightarrow E_B=\End_B(P_B)
\]
is a dg quasi-isomorphism.

Take the auxiliary grading group in \cite[Proposition 4.5(1)]{Fushimi} to be trivial. For $C=A,B$, the hypotheses imply
\[
 \perf C\subset\calD^-_{\mathrm{lf}}C
 :=\{X:H^{\gg0}X=0,\ \dim_kH^iX<\infty\ \forall i\},
\]
since this inclusion holds for $C$ and is preserved by finite cones, shifts and retracts. The simple sum in that proposition is $l$. Its contravariant Koszul functor
$K_C=\RHom_C(-,l)$ is therefore fully faithful on $\perf C$.
Under $E_A\simeq E_B$, derived adjunction identifies
\[
 K_BF(X)=\RHom_B(X\Lox_AB,l)\simeq\RHom_A(X,l)=K_A(X).
\]
This is an identification of Koszul-dual modules: the action is induced by the same extension map $E_A\to E_B$. Consequently, for $X,Y\in\perf A$,
\[
 \begin{aligned}
 \Hom_{\calD B}(FX,FY[n])
 &\iso\Hom(K_B(FY[n]),K_B(FX))\\
 &\iso\Hom(K_A(Y[n]),K_A(X))
 \iso\Hom_{\calD A}(X,Y[n]).
 \end{aligned}
\]
The composite is inverse to extension on Hom. Taking $X=Y=A$ identifies the canonical map
$H^n(A)\to H^n(B)$ with an isomorphism for every $n$.
\end{proof}

\begin{prp}\label{prop:canonical-completion}
Let $\La$ be basic, finite dimensional, $\ta_2$-finite and of global dimension at most two. Write $\La=\bbC Q/I$ with admissible $I$ and choose minimal relations $r$. Let $\widetilde Q$ be obtained by adding the relation arrows $\ro_r$, and put $W=\sum_r\ro_r r$. Then there exists a zigzag of quasi-isomorphisms of ordinary dg algebras
\[
 \Pi_3(\La)\simeq
 \Ga_{\mathrm{alg}}(\widetilde Q,W)
 \longrightarrow\widehat\Ga(\widetilde Q,W),
\]
where the last arrow is the canonical completion map. Its derived extension sends the regular module to the regular module.
\end{prp}
\begin{proof}
The first quasi-isomorphism follows from Keller's description of the undeformed completion \cite[Theorem 6.10]{Keller}. The ordinary Ginzburg algebra $G$ is quasi-free on the finite graded arrow space in degrees $0,-1,-2$. The relations have length at least two, so $W$ has length at least three and its differential is decomposable.

Its $H^0$ is the finite-dimensional algebra
$T_\La\operatorname{Ext}^2_\La(D\La,\La)$.
The positive tensor-degree ideal is nilpotent, and the quotient by it is $\La$. Thus the ideal generated by the old and relation arrows is the radical of $H^0G$, with quotient $l=\bbC^{Q_0}$. In particular the ordinary Jacobian ideal contains a power of the arrow ideal. Passing to completion does not change this quotient, so $H^0\widehat G\iso H^0G$ and has the same semisimple quotient.

The ordinary algebra $G$ is smooth and connective with finite-dimensional $H^0$, hence cohomologically locally finite. The same local finiteness holds for the Jacobi-finite completed Ginzburg algebra; see the proof of \cite[Theorem 4.17]{HK}. \lemref{lem:completion} applies. The statement about the regular module follows from extension of scalars.
\end{proof}

\begin{rmk}
For a tensor model $\Pi=T_\La\Th$, its tensor-length completion is
$\widehat\Pi_{\mathrm{ten}}=\prod_{m\ge0}\Th^{\ox_\La m}$ degree by degree. The differential preserves tensor length. Exactness of direct sums and products over $k$, and local finiteness, give
\[
 H^n\Pi=\Ds_{m\ge0}H^n(\Th^{\ox m})
 \ya{\sim}\prod_{m\ge0}H^n(\Th^{\ox m})
 =H^n\widehat\Pi_{\mathrm{ten}}.
\]
Indeed only finitely many summands are nonzero for each $n$. Thus $\Pi\to\widehat\Pi_{\mathrm{ten}}$ is a quasi-isomorphism. By \prpref{prop:canonical-completion}, we therefore obtain a zigzag of quasi-isomorphisms to the arrow completion. This does not identify their topologies or their underlying algebras.
\end{rmk}

\subsection{Pseudocompact Morita data for a specified silting object}
\begin{lem}\label{lem:pc-perfect}
For a fixed nonpositive augmented pseudocompact dg algebra $C\in\mathrm{PCAlg}_c(l)$ over a finite product $l$ of copies of $\bbC$, forgetting topology identifies its perfect dg category with the ordinary perfect dg category of its underlying dg algebra. The comparison preserves the regular object, finite-cell objects and their derived endomorphism dg algebras.
\end{lem}
\begin{proof}
For finite-cell modules $P,Q$, homogeneous maps are finite matrices
with entries in $C$, hence continuous. Therefore
\[
 \Hom^\bullet_{C,\mathrm{cont}}(P,Q)
 =\Hom^\bullet_{\Fgt_{\mathrm{top}} C}(\Fgt_{\mathrm{top}} P,\Fgt_{\mathrm{top}} Q).
\]
Finite-cell modules are cofibrant in both settings, and every
pseudocompact module is fibrant. The equality thus computes the
derived Hom complexes. The comparison sends $C$ to $C$ and extends
through
\[
 \{C\}\ \longmapsto\
 \{\text{finite sums, shifts and cones of }C\}
 \ \longmapsto\ \operatorname{per}C,
\]
where the last operation is dg idempotent completion. This proves
the assertion for perfect objects.
\end{proof}

\begin{lem}\label{lem:pointed-pc}
Let $C$ be a completed Jacobi-finite Ginzburg algebra over $\bbC$, with a reduced potential, and let $M$ be a basic silting object of $\perf C$. Assume either that $M$ is two-term, or that a finite-cell representative $P$ is given whose graded module is a finite sum of shifts of vertex projectives and whose twisting differential has entries in the arrow ideal. Then there are a nonpositive augmented pseudocompact dg algebra $A\in\mathrm{PCAlg}_c(l')$, for a finite product $l'$ of copies of $\bbC$, and a continuous $A$--$C$ dg bimodule $P$ such that:
\begin{enumerate}[label=(\roman*)]
\item $P_C$ represents the specified object $M$;
\item the underlying ordinary dg algebra of $A$ is quasi-isomorphic to $\REnd_C(M)$;
\item completed derived tensor with $P$ is a pseudocompact derived equivalence and sends the regular object $A$ to $M$.
\end{enumerate}
\end{lem}
\begin{proof}
Write $\calD_{\mathrm{pc}}(R)$ for the pseudocompact derived category
of \cite[Section 6]{VdB}. We construct $A$, $P$ and
$Q=\Hom^\bullet_{C,\mathrm{cont}}(P,C)$, and prove the bimodule weak
equivalences
\begin{equation}\label{eq:pointed-morita-target}
 A\longrightarrow P\widehat\ox_C^{\mathbf L}Q,
 \qquad Q\widehat\ox_A^{\mathbf L}P\longrightarrow C.
\end{equation}
For a nonpositive augmented pseudocompact dg algebra $R$, we use
the dual of \cite[Lemma A.2.1]{VdB} in the form
\begin{equation}\label{eq:pc-bounded-test}
 \left.
 \begin{gathered}
 X^{>u}=Y^{>v}=0,\\
 H^*(f):H^*X\xrightarrow{\sim}H^*Y
 \end{gathered}
 \right\}
 \quad\Longrightarrow\quad
 f:X\to Y\text{ is a weak equivalence}.
\end{equation}
The bounds will be checked separately for the two maps in
\eqref{eq:pointed-morita-target}.

\medskip\noindent\emph{Step 1.}
In the two-term case, choose a triangle $P_1\ya{f}P_0\to M\to P_1[1]$, where $P_i$ are finite sums of vertex projectives. Indeed, $H^0$ identifies $\add C$ with the finitely generated projectives over the basic algebra $H^0C$. A degree-zero morphism between these representatives is a matrix of elements of $C^0$; every such element is a cycle because $C^1=0$.

Let $\mathfrak m$ be the ideal of all graded arrows of $C$. Still in the two-term case, if the reduction of $f$ modulo $\mathfrak m$ is nonzero, elementary changes of basis isolate an invertible block between copies of the same vertex projective. Its inverse exists in the completed corner algebra by a geometric series. Row and column elimination splits off a contractible cone. Repeating finitely many times produces a representative with $f\equiv0\pmod{\mathfrak m}$. Thus cancelling these blocks changes neither the object $M$ nor its derived endomorphism algebra. Set $P=\Cone(f)$ in this case; otherwise use the specified minimal finite-cell representative. Put
\[
 E=\End^\bullet_{C,\mathrm{cont}}(P),\qquad
 A=\ta_{\le0}E.
\]
Here the truncation is the dg subalgebra with $A^i=E^i$ for $i<0$, $A^0=Z^0E$ and $A^i=0$ for $i>0$. Let $[a,b]$ contain the degrees of the finitely many generators of $P$. Then
\[
 P^{>b}=0,\qquad E^{>b-a}=0.
\]
For the two-term representative one may take $a=-1$, $b=0$.

Since $P$ is a finite-cell representative, it computes derived Hom. By expressing its endomorphisms as finite matrices, we obtain an equality of dg algebras
\[
 \End^\bullet_{C,\mathrm{cont}}(P)
 =\End^\bullet_{\Fgt_{\mathrm{top}}(C)}(\Fgt_{\mathrm{top}}(P)).
\]
Indeed every matrix entry acts continuously, and only finitely many entries occur. Since $M$ is silting,
\[
 H^iE=\Hom_{\perf C}(M,M[i])=0\quad(i>0).
\]
Consequently
\[
 A\hookrightarrow E\simeq\REnd_C(M),\qquad P_C\simeq M,
\]
which proves (i) and (ii).

\medskip\noindent\emph{Step 2.}
Let $F^nE$ be the subspace of matrices with entries in $\mathfrak m^n$. Since the potential is reduced, $d_C$ takes every arrow to a sum of paths of length at least two. Since the twisting matrix of $P$ has entries in $\mathfrak m$, the matrix differential satisfies
\[
 d_E(F^nE)\subset F^{n+1}E,\qquad
 F^aE\,F^bE\subset F^{a+b}E.
\]
In particular every $F^nE$ is a dg ideal. The quotients are finite dimensional as total graded vector spaces: they involve only finitely many paths of length less than $n$ and finitely many matrix positions. The matrix algebra is complete and separated. Its subspace $Z^0E$ is closed, so $A$ is pseudocompact for the induced filtration $F^nA=A\cap F^nE$.

Set $J_0=d(A^{-1})$ and $I=F^1A^0$. Then $J_0\subset I$. For $x\in I$ and $y\in A^0$,
\[
 (yx)^j\in F^jA^0,\qquad
 (1-yx)^{-1}=\sum_{j\ge0}(yx)^j\in A^0.
\]
By completeness, the series converges. Hence $I\subseteq\rad A^0$. The image $J_0$ is closed because it is the image of a continuous linear map between linearly compact spaces. In particular $H^0A=A^0/J_0$ has its quotient topology. It is finite dimensional: $C$ has finite-dimensional cohomology in each degree, and taking finite cones and retracts preserves this property for the Hom complexes of perfect objects. Since $J_0\subset\rad A^0$, it follows from the radical formula for a quotient that
\[
 A^0/\rad A^0\iso
 H^0A/\rad H^0A=:l'.
\]
Since $M$ is basic and $\perf C$ is Krull--Schmidt over the algebraically closed field $\bbC$, we have $l'\iso\bbC^r$.

To construct a section, put $A_n=A^0/F^nA^0$. Then
\[
 \begin{gathered}
 A_n/\rad A_n\cong l',\qquad \ker(A_{n+1}\to A_n)^2=0,\\
 \ker(A_{n+1}\to A_n)=F^nA^0/F^{n+1}A^0.
 \end{gathered}
\]
Lift the standard idempotents of $l'$ first to $A_1$, and then
successively along $A_{n+1}\to A_n$. The resulting systems satisfy
\[
 e_{i,n}e_{j,n}=\delta_{ij}e_{i,n},\qquad
 \sum_i e_{i,n}=1,\qquad e_{i,n+1}\longmapsto e_{i,n}.
\]
Hence $e_i=\varprojlim_n e_{i,n}$ defines a section $l'\to A^0$.

The map $A\to H^0A\to l'$ is continuous: $J_0$ is closed and the finite-dimensional Hausdorff quotient $H^0A$ is discrete. It is a dg augmentation split by the preceding section. In every finite-dimensional graded quotient $A/F^nA$, the ideal of negative degrees is nilpotent. Its quotient is $A^0/F^nA^0$, whose semisimple quotient is $l'$. Thus the kernel of the induced augmentation is exactly the radical of each such finite quotient. Passing to the inverse limit identifies $\ker(A\to l')$ with the topological graded radical. Thus $A\in\mathrm{PCAlg}_c(l')$.

\medskip\noindent\emph{Step 3.}

\smallskip\noindent\emph{(a) The adjunction and its counit at $C$.}
Restrict the left action of $E$ on $P$ to $A$, and put
\[
 Q=\Hom^\bullet_{C,\mathrm{cont}}(P,C).
\]
Then $P$ is an $A$--$C$ bimodule and $Q$ is a $C$--$A$ bimodule. All actions are continuous matrix operations. The generators of the left $C$-module $Q$ have degrees in $[-b,-a]$, so $Q^{>-a}=0$. By finite-projective duality, we obtain natural isomorphisms of complexes
\[
 Y\widehat\ox_C Q
 \iso\Hom^\bullet_{C,\mathrm{cont}}(P,Y),\qquad
 P\widehat\ox_C Q\iso E.
\]
These isomorphisms commute with the additional module actions and the differentials. The first is the evaluation isomorphism for a shifted vertex projective and extends to $P$ by finite direct sums and cones.

Write $\mathcal N_C,\mathcal N_A$ for the null classes of the pseudocompact model structures in \cite[Section 6]{VdB}. Since $P$ is finitely generated projective as a graded module, we have
\[
 \begin{gathered}
 G_0=\Hom^\bullet_{C,\mathrm{cont}}(P,-)\text{ exact on underlying graded modules},\\
 G_0(\prod_jX_j)\iso\prod_jG_0(X_j),\qquad
 G_0(\Cone f)\iso\Cone(G_0f),\\ G_0(X[1])\iso G_0(X)[1].
 \end{gathered}
\]
It preserves retracts and total complexes of short exact sequences. Since these total complexes generate the null class under the displayed operations,
\[
 G_0(\mathcal N_C)\subset\mathcal N_A.
\]
For a trivial fibration $f:X\twoheadrightarrow Y$, it follows from exactness that
\[
 G_0(f)\text{ surjective},\qquad
 \ker G_0(f)=G_0(\ker f)\in\mathcal N_A.
\]
Hence $G_0$ preserves fibrations and trivial fibrations. The completed tensor--Hom adjunction is Quillen, with derived functors
\[
 F=-\widehat\ox_A^{\mathbf L}P,\qquad G=G_0.
\]
Every object is fibrant, so no replacement is needed for $G$. This step applies without a boundedness assumption.

The unit at the cofibrant regular module $A$ is precisely $A\to G(P)=E$. It is a quasi-isomorphism by Step 1. Since $A^{>0}=0$ and $E^{>b-a}=0$,
\eqref{eq:pc-bounded-test} makes it a weak equivalence of right
$A$-modules. By the adjunction identity, we have
\[
 \ep_P F(\et_A)=1_P,
 \qquad \ep_P=F(\et_A)^{-1}.
\]
The full subcategory $\mathcal U=\{Y:\ep_Y\text{ is invertible}\}$ is thick. Therefore
\[
 \begin{gathered}
 P\in\mathcal U\ \Longrightarrow\ C\in\thick(P)\subseteq\mathcal U,\\
 Q=G(C)\in\thick(G(P))=\thick(A)=\perf A.
 \end{gathered}
\]

\medskip\noindent\emph{(b) A bimodule resolution and the same evaluation map.}
Consider the following maps induced by the counit and the action, respectively:
\[
 \ep_C:Q\widehat\ox_A^{\mathbf L}P\longrightarrow C,
 \qquad
 \et_A:A\longrightarrow P\widehat\ox_C^{\mathbf L}Q.
\]
The first is a morphism of $C$--$C$ bimodules and the second of $A$--$A$ bimodules. For the second map, the derived tensor product is represented by $P\widehat\ox_C Q=E$, since $P$ is finite-cell over $C$; the comparison is the inclusion $A\to E$.

To compute the first map in the bimodule category, we use the following resolution. Let $R\in\mathrm{PCAlg}_c(S)$ be nonpositive, where $S$ is a finite product of copies of $\bbC$, and let $X^{>b}=0$. A completed normalized bar resolution has underlying graded module
\[
 \operatorname{Bar}_R(X)
 =\prod_{n\geq0}
 X\widehat\ox_S (\overline R[1])^{\widehat\ox_S n}
       \widehat\ox_S R,
 \qquad \overline R=\ker(R\to S).
\]
The differential is the sum of the internal differentials and the bar differential induced by multiplication and the action on $X$. The augmentation to $X$ is induced by the module action on the length-zero term. The usual extra degeneracy defines an $S$-linear contraction of its augmented bar complex; normalization preserves this quasi-isomorphism. All these operations are continuous. Put $B_n=X\widehat\ox_S(\overline R[1])^{\widehat\ox_S n}\widehat\ox_S R$. Since $R^{>0}=0$,
\[
 B_n^j\ne0\ \Longrightarrow\ j\leq b-n,
 \qquad
 \operatorname{Bar}_R(X)^j=\prod_{0\leq n\leq b-j}B_n^j,
 \qquad \operatorname{Bar}_R(X)^{>b}=0.
\]
Thus each degree involves finitely many bar lengths.

Its underlying graded module is topologically free relative to $S$, hence graded projective. Indeed continuous duality makes the category of pseudocompact $S$-modules the opposite of the category of discrete $S$-modules, which is semisimple. Induction $V\mapsto V\widehat\ox_S R$ takes such projectives to graded projectives by its adjunction with restriction. Thus the bar module is cofibrant by the model structure in \cite[Section 6]{VdB}. Its bounded-above augmentation is a weak equivalence by the boundedness criterion, not by an assertion about arbitrary unbounded acyclic modules.

Apply the bar construction over
\[
 D=C^{\mathrm{op}}\widehat\ox_k A,
 \qquad S=l^{\mathrm{op}}\ox_k l'.
\]
Right $D$-modules are $C$--$A$ bimodules. This is a nonpositive
augmented pseudocompact dg algebra. Indeed, for finite-dimensional
quotients $C_i,A_j$, put
\[
 I_{ij}=\rad(C_i)^{\mathrm{op}}\ox A_j
       +C_i^{\mathrm{op}}\ox\rad(A_j).
\]
If the two radical nilpotence indices are $u_i,v_j$, then
\[
 I_{ij}^{u_i+v_j}=0,\qquad
 (C_i^{\mathrm{op}}\ox A_j)/I_{ij}\cong S.
\]
By passing to the inverse limit, we obtain the required topological radical
condition on $D$.

Let
\[
 q:\widetilde Q=\operatorname{Bar}_D(Q)\longrightarrow Q.
\]
Since $Q^{>-a}=0$, it follows from the preceding bar-degree estimate that
$\widetilde Q^{>-a}=0$. Its underlying graded module is induced
from $S$; after restriction to $A$ its induced terms have the form
\[
 V\widehat\ox_S D
 \cong W\widehat\ox_{l'} A
\]
for a pseudocompact graded $l'$-module $W$ (the $C$-factor is
absorbed into $W$). Since $l'$ is semisimple, induction and its
restriction adjunction give
\[
 \operatorname{Res}_A\widetilde Q
 \text{ graded projective}
 \quad\Longrightarrow\quad
 \operatorname{Res}_A\widetilde Q\text{ cofibrant}.
\]
Moreover, $q$ is an ordinary quasi-isomorphism of bounded-above
complexes. Applying \eqref{eq:pc-bounded-test} both over $D$ and
over $A$ shows that $q$ and $\operatorname{Res}_Aq$ are weak
equivalences. Thus the same resolution computes the tensor product
in bimodules and in right modules.

Consider the continuous bimodule map
\[
 \widetilde\varepsilon:
 \widetilde Q\widehat\ox_A P
 \xrightarrow{q\widehat\ox 1}Q\widehat\ox_A P
 \xrightarrow{\mathrm{ev}}C.
\]
By the restriction calculation, its underlying right $C$-module
map represents the counit $\varepsilon_C$ proved invertible in (a).
Consequently
\[
 H^*\operatorname{Cone}(\widetilde\varepsilon)=0,
 \qquad
 (\widetilde Q\widehat\ox_A P)^{>b-a}=0,
 \qquad C^{>0}=0.
\]
The source bound follows from
$\widetilde Q^{>-a}=0$ and $P^{>b}=0$. The cone is therefore
bounded above. The same finite-quotient calculation as for $D$
applies to $C^{\mathrm{op}}\widehat\ox_k C$, so
\eqref{eq:pc-bounded-test} makes $\widetilde\varepsilon$ a weak
equivalence of $C$--$C$ bimodules.

For the other comparison, since $P$ is finite-cell over $C$, we obtain
\[
 P\widehat\ox_C^{\mathbf L}Q
 \simeq P\widehat\ox_CQ\cong E.
\]
The action map $\eta_A:A\hookrightarrow E$ is an ordinary
quasi-isomorphism of $A$--$A$ bimodules, with
$A^{>0}=0$ and $E^{>b-a}=0$. Apply
\eqref{eq:pc-bounded-test} over
$A^{\mathrm{op}}\widehat\ox_k A$. This proves both bimodule
weak equivalences in \eqref{eq:pointed-morita-target}.

\medskip\noindent\emph{Step 4.}
Derived completed tensor is associative, with the regular bimodule as unit. The two bimodule weak equivalences give, for all objects $X$ and $Y$ of the respective pseudocompact derived categories,
\begin{align*}
 (X\widehat\ox_A^{\mathbf L}P)
       \widehat\ox_C^{\mathbf L}Q
 &\iso X\widehat\ox_A^{\mathbf L}A\iso X,\\
 (Y\widehat\ox_C^{\mathbf L}Q)
       \widehat\ox_A^{\mathbf L}P
 &\iso Y\widehat\ox_C^{\mathbf L}C\iso Y.
\end{align*}
Consequently completed tensor with $Q$ is inverse to $F$, and $F(A)=P$ represents the original $M$. This proves (iii). The bimodule identities establish the equivalence on the full pseudocompact derived categories. The generation of $\perf C$ was used to prove that the counit at $C$ is invertible.
\end{proof}

\begin{rmk}
In the second case, existence of the specified minimal finite-cell representative is part of the hypothesis. A perfect object is initially a retract of a finite-cell object. The lemma does not identify its endomorphism algebra with that of its cluster image outside the fundamental domain.
\end{rmk}

\subsection{Potential classes and ordinary dg recognition}
\begin{lem}\label{lem:cyclic-zero}
Let $D=\widehat\Ga(Q,W)$ be a nonpositive Jacobi-finite completed Ginzburg algebra, augmented over $l=\bbC^{Q_0}$, and suppose $W$ has no constant term. The class of $W$ in reduced continuous $HC_0(D)$ is zero if and only if its class in $HH_0(H^0D)$ is zero.
\end{lem}
\begin{proof}
Set $J=H^0D$. In the homological total grading of the continuous cyclic complex, an element in Hochschild length $n$ and internal cohomological degree $-j$ has total degree $n+j$. Both $n$ and $j$ are nonnegative. Thus total degree zero comes only from $D^0$ in length zero. Boundaries into it are $d(D^{-1})$ and commutators of degree-zero elements; higher cyclic columns have strictly positive total degree. By the definition of reduced continuous cyclic homology, we obtain
\[
 HC^{\mathrm{red}}_0(D)
 =D^0\big/\bigl(d(D^{-1})+
       \overline{[D^0,D^0]}+l\bigr)
 \iso J/([J,J]+l).
\]
Here commutators denote their linear span, not the two-sided ideal they generate. The image of $d$ is closed. Since $J$ is finite dimensional, taking the closure of the commutator space does not change its image in $J$. This also explains the degree-zero formula used in \cite[Section 4.5]{HK}.

Suppose first that the reduced cyclic class is zero. For the unit inclusion $\io:l\to HH_0(J)$ and augmentation retraction $\epsilon_*$,
\[
 [W]\in\operatorname{im}\io,\qquad
 \epsilon_*\io=1_l,\qquad \epsilon_*[W]=0.
\]
Consequently $[W]=\io\epsilon_*[W]=0$ in $HH_0(J)$. The converse follows from the quotient map $J/[J,J]	o J/([J,J]+l)$.
\end{proof}

\begin{prp}\label{prop:zero-class-recognition}
Let $C_i=\widehat\Ga(Q_i,W_i)$ be reduced Jacobi-finite completed Ginzburg algebras over $\bbC$, put $J_i=H^0C_i$, and assume
\[
 [W_i]=0\quad\text{in }HH_0(J_i)\qquad(i=1,2).
\]
If $M\in\perf C_1$ is basic two-term silting and
$\End_{\calC_{C_1}}(\pi M)\iso J_2$, then its ordinary derived endomorphism dg algebra is quasi-isomorphic to the underlying ordinary dg algebra of $C_2$.
Here $\calC_{C_i}$ is the topological generalized cluster category of \cite[Section 4]{HK}, with perfect objects identified by \lemref{lem:pc-perfect}.
\end{prp}
\begin{proof}
\noindent\emph{Step 1.}
Put $C=C_1$ and $B=\REnd_{C_1}(M)$. By \lemref{lem:pointed-pc}, there exist a nonpositive pointed pseudocompact algebra $A$ and a continuous Morita bimodule $P$ such that $P$ represents $M$ and $B\simeq\Fgt_{\mathrm{top}}(A)$ as ordinary dg algebras.

Let
\[
 \mathcal F=(\add C)*(\add C)[1].
\]
The restriction $\pi|_{\mathcal F}$ to the fundamental domain is an additive equivalence; see \cite[Proposition 4.2]{KY} and the fundamental-domain formulation recalled before \cite[Definition 5.1]{HI}. \lemref{lem:pc-perfect} identifies the morphisms between these perfect objects with their ordinary derived morphisms. Since $P$ represents the two-term object $M$, we have $P\in\mathcal F$. By full faithfulness of $\pi|_{\mathcal F}$ and the hypothesis on $\End_{\calC_{C_1}}(\pi M)$, we obtain
\[
 H^0A\iso\End_{\perf C}(P)
 \iso\End_{\calC_{C_1}}(\pi M)
 \iso J_2.
\]
The first isomorphism is induced by $A\to\End_C(P)$, the second by $\pi$, and the last is the isomorphism in the hypothesis.

\medskip\noindent\emph{Step 2.}
\lemref{lem:cyclic-zero} identifies reduced continuous degree-zero cyclic homology with $J_i/([J_i,J_i]+l_i)$. Hence the assumed equalities $[W_i]=0$ in $HH_0(J_i)$ imply
\[
 [W_i]=0\quad\text{in }HC^{\mathrm{red}}_0(C_i).
\]
We apply \cite[Corollary 4.27]{HK} to $C$, $A$ and $P$. By \lemref{lem:pointed-pc},
\[
 A\in\mathrm{PCAlg}_c(l'),\qquad A^{>0}=0,\qquad
 -\widehat\ox_A^{\mathbf L}P:
 \calD_{\mathrm{pc}}(A)\ya{\sim}\calD_{\mathrm{pc}}(C).
\]
Thus the hypotheses of \cite[Corollary 4.27]{HK} are satisfied, and there exists a weak equivalence
$s:\widehat\Ga(Q',W')\to A$ for which the cyclic homology isomorphism associated with $sP$ takes $[W']$ to $[W_1]$. This cyclic comparison is induced by the Morita bimodule $P$, as specified in \cite[Corollary 4.27]{HK}. That corollary transports the nondegenerate degree-two cyclic class and uses naturality of the periodicity map $S$ to identify the degree-zero potential classes. In particular, the induced map is an isomorphism
\[
 HC^{\mathrm{red}}_0(sP):
 HC^{\mathrm{red}}_0(\widehat\Ga(Q',W'))
 \ya{\sim} HC^{\mathrm{red}}_0(C),
 \qquad [W']\longmapsto[W_1]=0.
\]
Thus
\[
 [W']\in\ker HC^{\mathrm{red}}_0(sP)=0
 \quad\Longrightarrow\quad
 [W']=0\in HH_0(H^0\widehat\Ga(Q',W'))
\]
by \lemref{lem:cyclic-zero}. This uses the zero potential class of $C_1$. The minimal-model construction in \cite[Theorem 11.2.1]{VdB}, used in this recognition, has decomposable differential. Its Ginzburg presentation can therefore be taken with potential of order at least three. This is the reduced presentation used in the next step. Its Jacobian algebra is finite dimensional because the weak equivalence $s$ induces an isomorphism on $H^0$.

\medskip\noindent\emph{Step 3.}
The induced isomorphism on $H^0$ preserves the potential classes, since both are zero. Let $\al:H^0\widehat\Ga(Q',W')\to J_2$ be the isomorphism just constructed. After permuting the vertices, its images of the vertex idempotents and the chosen target vertex idempotents are two complete primitive orthogonal systems with the same simple labels. They are conjugate by a unit of the finite-dimensional target algebra. Thus, after a vertex permutation, there is a unit $u$ such that
\[
 \al'=\operatorname{Ad}(u)\al,\qquad
 \al'(e_i)=e_i,\qquad
 HH_0(\al')([W'])=[W_2]=0.
\]
The last equality uses $HH_0(\operatorname{Ad}(u))=1$.

Since both potentials have order at least three, their Jacobian relations belong to the square of the arrow ideal. Consequently the arrow spaces are the vertex components of the radical modulo its square. The marked isomorphism $\al'$ therefore identifies the two finite quivers after choosing bases of those spaces. Over $\bbC$, the reduced potentials and the marked isomorphism of their finite-dimensional Jacobian algebras therefore satisfy the hypotheses of noncommutative Mather--Yau. It follows that there is a right equivalence with the target potential and an isomorphism of the completed Ginzburg models \cite[Theorems 3.5 and 4.3]{HZ}. We next regard these maps as maps of ordinary dg algebras. The map $A\to\End_C(P)$ is the quasi-isomorphism of \lemref{lem:pointed-pc}, and forgetting topology identifies the latter with the ordinary derived endomorphism algebra of $M$. The map $s$ is a weak equivalence between nonpositive pseudocompact dg algebras. In this bounded-above setting it is an ordinary quasi-isomorphism on forgetting topology \cite[Section 6 and Appendix A]{VdB}. The right equivalence induces an actual continuous dg isomorphism between the two Ginzburg algebras. We therefore obtain the following zigzag of quasi-isomorphisms of ordinary dg algebras
\[
 B\simeq\Fgt_{\mathrm{top}}(A)
 \simeq\widehat\Ga(Q',W')
 \iso C_2.
\]

The quasi-isomorphism type on the left is that of the specified object $M$.
\end{proof}

\medskip\noindent\textbf{Ordinary recognition and prescribed gradings.}
\prpref{prop:zero-class-recognition} identifies ordinary dg quasi-isomorphism types. It does not assert that the right equivalence preserves a prescribed internal grading. Hua--Zhou give right-equivalent Jacobi-finite rigid potentials with different weights and Jacobian Hilbert series for the same arrow-weight system \cite[Example 7.6]{HZFinite}. Their Theorem 7.2 assumes weight preservation on
\[
 \mathfrak m/(\mathfrak m^2+J(W'))
 \longrightarrow\mathfrak m/(\mathfrak m^2+J(W_2))
\]
and obtains a bigraded isomorphism of the radical-associated graded Jacobian algebras, together with equality of the stated graded invariants. Even under this hypothesis, the original Jacobian isomorphism need not be homogeneous \cite[Example 7.3]{HZFinite}. Accordingly, the ordinary comparison proved above is used only as input to \thmref{thm:common-dg}, which separately constructs the two graded dg comparisons.

\begin{cor}\label{prop:recognition}
Let $\La_1,\La_2$ be basic finite-dimensional $\ta_2$-finite complex algebras of global dimension at most two, and put $\Ga_i=\Pi_3(\La_i)$. If $M\in\perf\Ga_1$ is basic two-term silting and
\[
 \End_{\calC_{\La_1}}(\pi M)\iso H^0\Ga_2,
\]
then $\REnd_{\Ga_1}(M)\simeq\Ga_2$ as ordinary dg algebras.
\end{cor}
\begin{proof}
\prpref{prop:canonical-completion} identifies $\Ga_i$ with the ordinary dg algebra underlying $C_i=\widehat\Ga(Q_i,W_i)$ for the canonical relation potentials. Extension of scalars preserves the regular objects, the specified $M$, its two-term property and its ordinary derived endomorphism algebra. For
\[
 W_i=\sum_r\ro_r r,
 \qquad d\Bigl(\sum_r\ro_r\ro_r^*\Bigr)=W_i,
\]
the potential class is zero. The fundamental-domain comparison identifies the endomorphism algebra of the transported object with $H^0\Ga_2=H^0C_2$. By \prpref{prop:zero-class-recognition}, we therefore obtain
\[
 \REnd_{\Ga_1}(M)\simeq\REnd_{C_1}(M\Lox_{\Ga_1}C_1)
 \simeq C_2\simeq\Ga_2.
\]
\end{proof}

\begin{rmk}
More generally, if $W$ is homogeneous of nonzero integer weight $d$ for arrow weights $w(a)$, then, by the cyclic Euler identity, we have
\[
 d[W]=\sum_a w(a)[a\partial_aW]=0
 \quad\text{in }HH_0(\Jac(Q,W)).
\]
Thus the zero-class hypothesis holds over $\bbC$.
\end{rmk}

\section{Covering quotients}\label{sec:cover}
In this section we identify the covering quotient of a bigraded dg algebra with a derived category of graded modules. We use one internal grading for the completion and the other for the grading on the degree-zero algebra.

For a bigraded dg algebra $R$, define
\[
 \mathcal Q(R)=\perf^{\bbZ^2}R/
 \bigl(\perf^{\bbZ^2}R\cap\calD_{\fd}^{\bbZ^2}R\bigr).
\]
The finite-dimensional condition includes all cohomological and internal degrees. The cover construction below uses ordinary direct-sum graded modules.

\subsection{Admissible models and their covering quotients}

Let $R$ be a dg algebra with internal gradings $(u,v)$.
Write $\mathcal R=\Cov_v(R)$ for its $v$-cover.
Thus
\[
 \operatorname{Ob}\mathcal R=\mathbb Z,
 \qquad
 \mathcal R(m,n)=R_{v=n-m}.
\]
The differential and composition are induced by those of $R$.
The $u$-grading induces an Adams grading on $\mathcal R$.
Right dg $\mathcal R$-modules identify with $v$-graded
right dg $R$-modules. We use direct-sum gradings throughout.

Suppose that $R$ is supported in $u\geq0$.
Let $\mathcal R_0$ be the dg category with the same objects
as $\mathcal R$ and morphism complexes
\[
 \mathcal R_0(m,n)=R_{u=0,v=n-m}.
\]
The projection $\mathcal R\to\mathcal R_0$ is a dg functor.
We regard the right and left representable modules
\[
 \mathcal R_0(-,n)
 \quad\text{and}\quad
 \mathcal R_0(n,-)
 \qquad(n\in\mathbb Z)
\]
as Adams-graded modules over $\mathcal R$ and
$\mathcal R^{\mathrm{op}}$, respectively, by restriction
along this projection.

We use the following terminology.
The Gorenstein parameter is understood in the sense of
\cite[Definition A.1]{Han}.

\begin{dfn*}\label{def:admissible}
Let $L$ be a finite-dimensional $\mathbb Z$-graded algebra.
A dg algebra $R$ with internal gradings $(u,v)$ is called
an \emph{admissible completion model} of $L$ if the following
conditions hold:
\begin{enumerate}[label=(\alph*)]
\item
$R^i=0$ for $i>0$, and $R_{u,v}=0$ for $u<0$.

\item
There is a zigzag of quasi-isomorphisms of $v$-graded
dg algebras
\[
 R_{u=0}\simeq L,
\]
where $L$ is concentrated in cohomological degree zero
and has finite global dimension.

\item
The Adams-graded dg category $\mathcal R=\Cov_v(R)$
has Gorenstein parameter one, and, for every $n\in\mathbb Z$,
\[
 \mathcal R_0(-,n)\in\perf^{\mathbb Z}_u\mathcal R,
 \qquad
 \mathcal R_0(n,-)\in
 \perf^{\mathbb Z}_u(\mathcal R^{\mathrm{op}}).
\]
\end{enumerate}
\end{dfn*}

By \cite[Lemma A.4]{Han}, the Gorenstein parameter condition
in {\rm(c)} is equivalent to
\[
 H^i\!\left(
 \underline{\operatorname{RHom}}_{\mathcal R}
 \bigl(\mathcal R_0(-,n),\mathcal R\bigr)
 \right)_u=0
 \qquad
 (i,n\in\mathbb Z,\ u\ne-1).
\]
Here the derived Hom is a left dg $\mathcal R$-module,
and the subscript denotes its Adams degree.
Together with the two perfectness conditions in {\rm(c)},
this is Setting~A.5 of \cite{Han} with parameter one.
The vanishing condition above characterizes the Gorenstein
parameter only. The two perfectness conditions in {\rm(c)}
are separate hypotheses. In the situation of
Theorem~\ref{thm:automatic}, they follow by transporting
the degree-zero modules along the singly graded
quasi-isomorphisms and then applying
Lemma~\ref{lem:forget}(i) on both sides.

\begin{prp}\label{prop:cover}
Let $R$ be an admissible completion model of $L$.
Then there is a triangle equivalence
\[
 \mathcal Q(R)\simeq\Db(\gr L).
\]
\end{prp}
\begin{proof}
\emph{Step 1.}
By Hanihara's Beilinson-type theorem with Gorenstein parameter one \cite[Theorem A.7]{Han}, we obtain
\begin{equation}\label{eq:han}
 \perf\mathcal R_0\simeq
 \perf^{\bbZ}_u\mathcal R/
 \thick\{\mathcal R_0(a)\mid a\in\bbZ\},
\end{equation}
where all representable modules of $\mathcal R_0$ are included. The numerator is $\perf^{\bbZ^2}R$. 

\medskip\noindent\emph{Step 2.}
Every representable Adams-zero module, viewed as a bigraded $R$-module, has finite-dimensional total cohomology: its cohomology is a shifted copy of the finite-dimensional algebra $L$. It is perfect by assumption. Consequently the denominator of \eqref{eq:han} is contained in $\perf^{\bbZ^2}R\cap\calD_{\fd}^{\bbZ^2}R$.

Conversely, let $X$ be perfect with finite-dimensional total cohomology. In the ambient derived category the standard cohomological truncation triangles express it using its finitely many nonzero cohomology modules over $H^0R$. We do not yet assume these truncations are perfect. The following filtration argument proves that each cohomology module belongs to the stated thick subcategory, and consequently proves its perfectness as well. For $M=H^iX$ put
\[
 \begin{gathered}
 F^aM=\Ds_{u\ge a,\,v\in\bbZ}M_{u,v},\qquad
 F^aM\,H^0R\subset F^aM,\\
 (F^aM/F^{a+1}M)(H^0R)_{u>0}=0.
 \end{gathered}
\]
Finite support makes the filtration finite. Every layer is therefore a finite-dimensional $v$-graded $L$-module, up to a $u$-shift.

Let $V$ be such a layer and put $d=\operatorname{gl.dim}L$. Successive finite graded free surjections give
\[
 0\longrightarrow K_d\longrightarrow F_{d-1}\longrightarrow\cdots
 \longrightarrow F_0\longrightarrow V\longrightarrow0,
 \qquad \Fgt(K_d)\in\prj L.
\]
By taking the degree-zero part of a splitting, we obtain $K_d\in\gr^v\prj L$ (for $d=0$, $V$ itself is graded projective). If $\mathcal T$ denotes the denominator of \eqref{eq:han}, we deduce from the finite filtrations that
\[
 V\in\mathcal T\ \Longrightarrow\ H^iX\in\mathcal T\ (\forall i)
 \ \Longrightarrow\ X\in\mathcal T.
\]
Thus $\mathcal T=\perf^{\bbZ^2}R\cap\calD_{\fd}^{\bbZ^2}R$.

\medskip\noindent\emph{Step 3.}
Since $\mathcal R_0\simeq\Cov_vL$ and $\operatorname{gl.dim}L<\infty$,
\[
 \mathcal Q(R)\simeq\perf\mathcal R_0
 \simeq\perf(\Cov_vL)
 \simeq\Kb(\gr\prj L)
 \simeq\Db(\gr L).
\]

\end{proof}

\subsection{Comparison of covering quotients}
\begin{thm}\label{thm:automatic}
Let $R$ be an ordinary bigraded connective dg algebra, with internal degrees $(u,v)$. Suppose there is a zigzag of quasi-isomorphisms of $u$-graded dg algebras
\[
 (R,u)\simeq(\Pi_3(\La),u),
\]
where $\La$ is finite dimensional of global dimension at most two and $u$ is the tensor grading. Then $R$ is bigraded quasi-isomorphic to an admissible completion model of a $v$-graded algebra $L$ whose underlying algebra is isomorphic to $\La$. In particular
\[
 \mathcal Q(R)\simeq\Db(\gr L).
\]
\end{thm}
\begin{proof}
\emph{Step 1.}
The cohomology in negative $u$-degrees vanishes. Consequently the subalgebra $R_{u\ge0}$ maps quasi-isomorphically to $R$. Replace it by its good nonpositive cohomological truncation. These replacements preserve both internal gradings. Apply the corresponding single-graded replacements to every intermediate algebra in the given zigzag. Restriction to $u$-degree zero induces a zigzag of quasi-isomorphisms
$R_0\simeq\La$ of ordinary dg algebras. Thus $R_0$ has cohomology only in degree zero, and its canonical map to $L=H^0R_0$ is a $v$-graded quasi-isomorphism. After forgetting $v$, we obtain $L\iso\La$.

\medskip\noindent\emph{Step 2.}
The ordinary tensor completion has its zero part perfect on both sides and has Gorenstein parameter one \cite[Example A.12(2)]{Han}. The natural augmentations to degree zero commute with nonnegative $u$-graded dg maps. For a single arrow $f:S\to T$ in this nonnegative zigzag, the map $S_0\to\operatorname{Res}_fT_0$ is a quasi-isomorphism of graded $S$-modules, because $f_0$ is a quasi-isomorphism and both actions factor through the zero-degree augmentations. Since extension along $f$ is inverse to restriction, it follows that
\[
 S_0\Lox_S T\simeq T_0.
\]
The same argument applies to left modules. Therefore extension along the given zigzag identifies the degree-zero modules with the degree-zero modules, not just with unspecified perfect objects. It follows that $R_0$ is left and right perfect after forgetting $v$. \lemref{lem:forget}(i), on each side, reflects this property back to the bigraded categories.

\medskip\noindent\emph{Step 3.}
Perfectness and \lemref{lem:forget}(ii) identify
\[
 \Fgt_v\uRHom_R(R_0,R)
 \simeq\underline{\RHom}_{\Fgt_vR}(\Fgt_vR_0,\Fgt_vR).
\]
For $u\ne-1$, the singly graded quasi-isomorphisms imply
\[
 0=H^i\!\left(\Fgt_v\uRHom_R(R_0,R)\right)_u
   =\Ds_v H^i\!\left(\uRHom_R(R_0,R)\right)_{u,v}.
\]
Hence every summand vanishes.

Let $\mathcal R$ be the $v$-cover. Its Adams-zero representables correspond to the internal $v$-shifts of $R_0$ and hence are perfect on both sides. The preceding vanishing implies that
$\RHom_{\mathcal R}(\mathcal R_0,\mathcal R)$ is supported in Adams degree $-1$. By \cite[Lemma A.4]{Han}, $\mathcal R$ has Gorenstein parameter one. Thus the conditions of Definition \ref{def:admissible} hold. The conclusion follows from \prpref{prop:cover}.
\end{proof}

\begin{prp}\label{prop:morita-descent}
Let $\widetilde N$ be a bigraded perfect $R$-module whose underlying object is a perfect generator, and let
$A=\uREnd_R^{\bbZ^2}(\widetilde N)$. There is a bigraded Morita equivalence induced by $\widetilde N$. This equivalence identifies the perfect subcategories of finite-dimensional total cohomology and hence induces
\[
 \overline G:\mathcal Q(A)\xrightarrow{\sim}\mathcal Q(R),
 \qquad
 \overline G\circ(a,b)\cong(a,b)\circ\overline G
 \quad((a,b)\in\bbZ^2).
\]
These isomorphisms are compatible with addition in $\bbZ^2$.
\end{prp}
\begin{proof}
The homogeneous Hom decomposition detects the zero object using the shifts of $\widetilde N$, so these shifts compactly generate. By applying the derived tensor--Hom criterion \cite[Lemma 3.10]{KellerDG} to the covering categories, we obtain a bigraded Morita equivalence. After forgetting the gradings, it is the Morita equivalence associated with the underlying perfect generator. To check the inverse functor, write $F=\uRHom_R(\widetilde N,-)$ and put $Q=F(R)$. The internal shifts of $Q$ are compact generators over $A$. Since $F$ is fully faithful, for every $X$ we have
\[
 \uRHom_A(Q,F(X))\simeq\uRHom_R(R,X)\simeq X
\]
is compatible with the right $R$-action. In particular, the inverse is also computed by derived Hom from a perfect object.

For $X\in\calD_{\fd}R$, the subcategory
\[
 \mathcal V_X=\{P\in\perf R:\ \dim_kH^*\RHom_R(P,X)<\infty\}
\]
is thick and contains $R$, hence equals $\perf R$. Apply this to both $\widetilde N$ and $Q$, using \lemref{lem:forget} to identify total homogeneous cohomology with its ungraded counterpart. Both functors preserve perfect objects with finite-dimensional total cohomology. They therefore descend to the stated Verdier quotients.
On a bigraded semi-free model, $\widetilde N$ is an $A$--$R$
dg bimodule, and the inverse of $F$ is
$G=-\Lox_A\widetilde N$. The canonical isomorphisms
\[
 G(X(a,b))\cong G(X)(a,b)
\]
are compatible with addition of internal degrees. The subcategories
being factored out are stable under both shifts, so these
isomorphisms descend to $\overline G$.
\end{proof}

\section{Proof of Theorem \texorpdfstring{\ref{thm:main}}{1.1}}\label{sec:main-proof}
We prove \thmref{thm:main} by applying the results of Sections~\ref{sec:relative}--\ref{sec:cover}.

Put $\Ga_i=\Pi_3(\La_i)$. The algebras $B,R,A$ constructed below enter the proof in the following order:
\[
\begin{gathered}
 \underbrace{B\simeq\Ga_2}_{\text{ordinary dg recognition: \corref{prop:recognition}}}
 \\
 \Downarrow\quad\text{\thmref{thm:common-dg}}
 \\
 \underbrace{(R,p)\simeq(B,p),\qquad(R,q)\simeq(\Ga_2,q)}_{
   \text{common bigraded model}}
 \\
 \Downarrow\quad\text{\thmref{thm:relative} and \corref{cor:replacement}}
 \\
 \underbrace{(A,p)\simeq(\Ga_1,p)}_{
   \text{relative lift retaining the first graded realization}}
 \\
 \Downarrow\quad\text{\prpref{prop:assembly}}
 \\
 \Db(\gr^qL_1)\simeq\mathcal Q(A)
       \simeq\mathcal Q(R)\simeq\Db(\gr^pL_2).
\end{gathered}
\]
The first comparison uses \lemref{lem:pointed-pc} and \prpref{prop:zero-class-recognition}. We shall also need the following consequence of Section~\ref{sec:cover}.

\begin{prp}\label{prop:assembly}
Let $\La_1,\La_2$ be as in \thmref{thm:main}. Suppose $R$ is an ordinary bigraded connective dg algebra and $\widetilde N$ is a bigraded perfect module whose underlying object is a perfect generator. Put $A=\uREnd_R^{\bbZ^2}(\widetilde N)$ and assume $A$ is connective. If there are zigzags of quasi-isomorphisms of singly graded dg algebras
\[
 (A,p)\simeq(\Pi_3(\La_1),p),\qquad
 (R,q)\simeq(\Pi_3(\La_2),q),
\]
where the gradings on the completions are their tensor gradings, then the two ordinary algebras admit gradings with equivalent bounded graded derived categories.
\end{prp}
\begin{proof}
\prpref{prop:morita-descent} identifies $\mathcal Q(A)$ and $\mathcal Q(R)$. Apply \thmref{thm:automatic} to $A$ with completion coordinate $p$ and residual coordinate $q$, and to $R$ with the coordinates reversed. We obtain
\[
 \Db(\gr^qL_1)\simeq\mathcal Q(A)\simeq\mathcal Q(R)
       \simeq\Db(\gr^pL_2),\qquad \Fgt L_i\iso\La_i.
\]
Transport the residual gradings through the last algebra isomorphisms.
\end{proof}

\begin{rmk}\label{rmk:residual-shifts}
The equivalence between the bigraded quotient categories
commutes with both internal shifts. However, the two
identifications in the proof use different completion
degrees:
\[
 \Db(\gr^{q}L_1)
 \simeq \mathcal Q(A)
 \simeq \mathcal Q(R)
 \simeq \Db(\gr^{p}L_2).
\]
The residual grading is $q$ on the left and $p$ on the right.
Thus compatibility with each of the shifts $(1,0)$ and
$(0,1)$ does not identify the grading shifts at the two
ends. The resulting triangle equivalence is not asserted
to commute with the internal grading shifts.
\end{rmk}

The objects and endomorphism models used below are related by
\[
 T\in\calC_{\La_1}
 \ \longleftarrow\ M\in\perf\Ga_1
 \ \longmapsto\ B=\uREnd^p_{\Ga_1}(\widetilde M)
 \ \longleftarrow\ R.
\]
We apply \thmref{thm:common-dg} to $B\simeq\Ga_2$ and denote the resulting model by $R$. The object to which \thmref{thm:relative} will be applied is the image of $\Ga_1[1]$ under the induced Morita equivalence.

\begin{proof}[Proof of \thmref{thm:main}]
\noindent\emph{Step 1.}
If one algebra is zero, its cluster category is zero. Since $\End(\pi\Ga_i)\iso H^0\Ga_i$ contains $\La_i$, cluster equivalence forces the other algebra to be zero as well, and the assertion is immediate. Hence we may assume both algebras are nonzero. We first prove the assertion for basic algebras. Put $\Ga_i=\Pi_3(\La_i)$, with canonical tensor gradings $p$ and $q$, respectively. They are connective, smooth and cohomologically locally finite. Let
\[
 E:\calC_{\La_2}\ya{\sim}\calC_{\La_1},
 \qquad T=E(\pi_2\Ga_2).
\]
For a connective $3$-Calabi--Yau dg algebra with finite-dimensional $H^0$, cluster-tilting objects lift to two-term silting objects in the fundamental domain; see \cite[Theorem 5.2]{HI}. Choose a basic two-term silting object $M\in\perf\Ga_1$ with $\pi_1M\iso T$. In particular
\[
 \End_{\calC_{\La_1}}(\pi_1M)\iso H^0\Ga_2.
\]

\medskip\noindent\emph{Step 2.}
Apply \thmref{thm:relative} with a trivial first grading to obtain a lift $\widetilde M$ for the canonical $p$-grading. Its internal shifts generate the graded derived category. Let
\[
 B=\uREnd^p_{\Ga_1}(\widetilde M).
\]
After forgetting $p$, this is the ordinary endomorphism dg algebra of $M$, by \lemref{lem:forget}. By \corref{prop:recognition}, $B$ and $\Ga_2$ are connected by a zigzag of quasi-isomorphisms of ordinary dg algebras. This comparison is taken after forgetting the internal gradings.

The ordinary algebra $B$ is connective because $M$ is silting. It is cohomologically locally finite because $M$ is perfect and $\Ga_1$ is cohomologically locally finite: the property of having finite-dimensional Hom in each degree extends from $\Ga_1$ to its thick closure in both arguments. Hence the hypotheses of \thmref{thm:common-dg} hold. By that theorem, there is a nonpositive bigraded algebra $R$ with
\begin{equation}\label{eq:common-ends}
 (R,p)\simeq(B,p),\qquad (R,q)\simeq(\Ga_2,q).
\end{equation}
\medskip\noindent\emph{Step 3.}
The $p$-graded Morita functor associated with $\widetilde M$ takes $\Ga_1[1]$ to a specified object of $\perf^pB$. Transport it along the first zigzag in \eqref{eq:common-ends} and denote the result by $N\in\perf^pR$.

For silting objects we use the order $P\ge Q$ if $\Hom(P,Q[i])=0$ for $i>0$. Since $M$ is two-term,
\[
 \Ga_1\ge M\ge\Ga_1[1].
\]
Let $\Ph:\perf\Ga_1\ya{\sim}\perf R$ be the underlying equivalence, with $\Ph(M)\iso R$ and $\Ph(\Ga_1[1])\iso N$. Then
\[
 \Ga_1\ge M\ge\Ga_1[1]
 \ \Longrightarrow\ M\ge\Ga_1[1]\ge M[1]
 \ \Longrightarrow\ \Ph(M)\ge\Ph(\Ga_1[1])\ge\Ph(M)[1].
\]
Thus $R\ge N\ge R[1]$, so the underlying $N$ is two-term silting relative to $R$. The finite-dimensional algebra $H^0R\iso H^0\Ga_2$ is basic. Apply \thmref{thm:relative} and \corref{cor:replacement} to construct a bigraded lift $\widetilde N$ with
\[
 \Fgt_q\widetilde N\simeq N\quad\text{in }\perf^pR.
\]
Set $A=\uREnd_R^{\bbZ^2}(\widetilde N)$. By the endomorphism comparison in \thmref{thm:relative}, we obtain quasi-isomorphisms preserving the $p$-grading:
\begin{equation}\label{eq:source-end}
 (A,p)\simeq
 \uREnd^p_{\Ga_1}(\Ga_1[1])
 \simeq(\Ga_1,p).
\end{equation}
The shift functor is a dg functor (with its usual signs on homogeneous morphisms), so it identifies these endomorphism dg algebras. In particular $A$ is connective, as required in \prpref{prop:assembly}.

\medskip\noindent\emph{Step 4.}
By \prpref{prop:morita-descent},
$\mathcal Q(A)\simeq\mathcal Q(R)$.
Apply \thmref{thm:automatic} to $A$ with $(u,v)=(p,q)$ and to $R$ with $(u,v)=(q,p)$, using \eqref{eq:source-end} and \eqref{eq:common-ends}. We obtain graded finite-dimensional algebras
\[
 L_1=H^0(A_{p=0}),\qquad L_2=H^0(R_{q=0})
\]
after the indicated quasi-isomorphic nonnegative replacements. Their underlying algebras are $\La_1,\La_2$, respectively. The residual gradings are $q$ on $L_1$ and $p$ on $L_2$. We obtain
\[
 \Db(\gr^qL_1)\simeq\mathcal Q(A)
 \simeq\mathcal Q(R)\simeq\Db(\gr^pL_2).
\]
Transport the gradings through the underlying algebra isomorphisms.

\medskip\noindent\emph{Step 5.}
First apply the construction to basic Morita representatives. Morita invariance of Calabi--Yau completion \cite[Proposition 4.2]{Keller} identifies their cluster categories with those of the original algebras, and preserves finite global dimension and $\ta_2$-finiteness. Each original finite-dimensional algebra is an endomorphism algebra of a finite direct sum of indecomposable projectives over its basic algebra. Equip each projective summand with the grading induced by a homogeneous primitive idempotent. The direct sum contains every indecomposable projective. It is therefore a graded projective generator when all its internal shifts are included. Its homogeneous endomorphism algebra acquires a grading, has the original algebra as its underlying algebra, and is graded Morita equivalent to the basic algebra. The assertion follows by composition with these graded Morita equivalences.
\end{proof}

\section{Internal shifts and the converse}\label{sec:scope}
Write $S_i=(1)$ for the internal shift on $\Db(\gr\La_i)$. Commutation with $S_i$ imposes a stronger condition than the conclusion of \thmref{thm:main}.

\begin{prp}\label{prop:shift-descent}
Let $\La_1,\La_2$ be finite-dimensional $\bbZ$-graded $k$-algebras of finite global dimension. Suppose there are a $k$-linear triangle equivalence
\[
 F:\Db(\gr\La_1)\ya{\sim}\Db(\gr\La_2)
\]
and a natural isomorphism of triangle functors $\et:FS_1\ya{\sim}S_2F$. Then the underlying algebras are derived equivalent. In particular, over $\bbC$, if they are $\ta_2$-finite and have global dimension at most two, they are cluster equivalent.
\end{prp}
\begin{proof}
By iterating $\et$ and its inverse, we obtain natural isomorphisms
\[
 \et_n:FS_1^n\ya{\sim}S_2^nF\quad(n\in\bbZ),
 \qquad \et_0=1,
 \qquad
 \et_{m+n,X}=S_2^m(\et_{n,X})\circ\et_{m,S_1^nX}.
\]
Since $\bbZ$ is generated freely by $1$, the iterated isomorphisms satisfy the unit and composition identities.

Put $T=F(\La_1)$ and $U=\Fgt(T)$. Since $\La_2$ has finite global dimension, $T$ has a bounded graded projective representative. Hence $U\in \Kb(\prj\La_2)$. For bounded graded projective complexes $X,Y$, by decomposing morphisms into homogeneous components, we obtain
\[
 \Hom_{\Db(\operatorname{mod}\La_2)}(\Fgt X,\Fgt Y[i])
 \iso\Ds_{n\in\bbZ}
 \Hom_{\Db(\gr\La_2)}(X,Y(n)[i]).
\]
The sums are finite on the representing Hom complexes. Applying this decomposition to $T$ and using $\et_n$, we obtain
\begin{align*}
 \Hom_{\Db(\operatorname{mod}\La_2)}(U,U[i])
 &\iso\Ds_n\Hom_{\Db(\gr\La_2)}(T,T(n)[i])\\
 &\iso\Ds_n\Hom_{\Db(\gr\La_1)}
              (\La_1,\La_1(n)[i])=0\qquad(i\ne0).
\end{align*}
For $i=0$ these identifications preserve composition. Indeed, a homogeneous map $f:\La_1\to\La_1(n)$ is sent to $\et_{n,\La_1}F(f)$; naturality and the coherence identity above preserve the product obtained by shifting one map before composition. Hence
\[
 \End_{\Db(\operatorname{mod}\La_2)}(U)
 \iso\Ds_n\Hom_{\Db(\gr\La_1)}
       (\La_1,\La_1(n))
 \iso\La_1
\]
as ordinary algebras. The last isomorphism is left multiplication on the right regular module.

Finally,
\[
 \Db(\gr\La_1)=\thick\{\La_1(n):n\in\bbZ\}
 \quad\Longrightarrow\quad
 \Db(\gr\La_2)=\thick\{T(n):n\in\bbZ\}.
\]
In particular, $\La_2$ lies in this thick closure and is obtained using finitely many such shifts. Since $\Fgt(T(n))\iso U$ for every $n$, after forgetting the grading we obtain
\[
 \La_2\in\thick(U),\qquad
 \thick(U)=\Kb(\prj\La_2)=\Db(\operatorname{mod}\La_2).
\]
Thus $U$ is a tilting complex with endomorphism algebra $\La_1$. By the tilting-complex theorem \cite[Theorem 6.4]{Rickard}, we obtain
\[
 \Db(\operatorname{mod}\La_1)
 \simeq\Db(\operatorname{mod}\La_2).
\]
Under the additional cluster-category hypotheses, it follows from \cite[Proposition 5.2]{AO} that $\calC_{\La_1}\simeq\calC_{\La_2}$. This conclusion asserts existence; no descent of the originally specified functor $F$ to a prescribed cluster functor is required.
\end{proof}

\begin{exm}\label{ex:no-converse}
Let $\La_1=\bbC$ and $\La_2=\bbC\times\bbC$, both trivially graded. They have global dimension zero and are $\ta_2$-finite. For a finitely generated graded vector space $V$, define
\[
 F(V)_m=(V_{2m},V_{2m+1})\qquad(m\in\bbZ),
\]
where the two factors carry the actions of the two primitive idempotents of $\La_2$. Conversely, for $W_m=(W_m^{(0)},W_m^{(1)})$, put
\[
 G(W)_{2m}=W_m^{(0)},\qquad
 G(W)_{2m+1}=W_m^{(1)}.
\]
These functors preserve finite support and are exact, with $GF=1$ and $FG=1$. By applying these functors termwise, we obtain
\[
 \gr\bbC\simeq\gr(\bbC\times\bbC),\qquad
 \Db(\gr\bbC)\simeq\Db(\gr(\bbC\times\bbC)).
\]
They do not commute with the internal shift. For example, if $V=\bbC$ is concentrated in degree zero, then $F(V(1))$ is supported in degree $-1$ at the second idempotent, whereas $F(V)(1)$ is supported there at the first idempotent. These two objects are not isomorphic.

For these semisimple algebras the Serre functor is the identity. Their classical and generalized cluster categories agree, and
\[
 \calC_{\bbC}\simeq\Db(\operatorname{mod}\bbC)/[2],
 \qquad
 \calC_{\bbC\times\bbC}\simeq\calC_{\bbC}\times\calC_{\bbC}.
\]
The first category has exactly two indecomposable isomorphism classes, represented by $\bbC$ and $\bbC[1]$, while the second has four. Thus
\[
 \Db(\gr\La_1)\simeq\Db(\gr\La_2)
 \qquad\text{but}\qquad
 \calC_{\La_1}\not\simeq\calC_{\La_2}.
\]
This example uses a disconnected algebra and makes no assertion about a converse restricted to connected algebras. It rules out the converse of \thmref{thm:main} in its stated generality.
\end{exm}

\section{Tilting lifts and acyclic cluster type}\label{sec:torus-lifting}
We consider cluster-tilting objects that admit tilting lifts whose endomorphism algebras have global dimension at most two. The results of this section are independent of \thmref{thm:main}. We use \cite[Theorem~3.5]{AO} and the surface results of \cite{AGSurface,AmiotTorus,BZ}. Equivalences and endomorphism algebras are ungraded unless otherwise stated.

For a $\ta_2$-finite algebra $A$ of global dimension at most two, put
\[
 \mathcal T_A=\left\{T\in\Db(A):
 \begin{array}{l}
 \Hom(T,T[n])=0\quad(n\ne0),\\
 \thick(T)=\Db(A),\quad \operatorname{gl.dim}\End(T)\le2
 \end{array}\right\}.
\]
The second part of \cite[Question 6.23, p.~47]{AOTree} asks whether
\[
 \pi_A(\mathcal T_A)=\CT(\calC_A)
 \quad\Longrightarrow\quad A\text{ is piecewise hereditary}.
\]
Here the equality concerns isomorphism classes; it is enough to consider basic objects. The published question imposes no acyclic-cluster-type hypothesis. We first prove the positive implication with that hypothesis, and then give a counterexample without it.

\subsection{The positive implication for acyclic cluster type}\label{subsec:acyclic-lifting}
For an algebra of acyclic cluster type, the cluster category contains a cluster-tilting object with hereditary endomorphism algebra. We show that a tilting lift of this object implies piecewise heredity.

\begin{lem}\label{lem:acyclic-canonical}
Let $B$ be a basic finite-dimensional $\ta_2$-finite $\bbC$-algebra of global dimension at most two. If the Gabriel quiver of
\[
 \widetilde B:=H^0\Pi_3(B)
       \iso\End_{\calC_B}(\pi_B B)
\]
is acyclic, then $B$ is hereditary and the canonical inclusion $B\to\widetilde B$ is an isomorphism.
\end{lem}
\begin{proof}
Choose an admissible presentation $B=\bbC Q_B/I$ and a minimal set $\mathcal R$ of relations, each with fixed source and target. The canonical relation-potential presentation, recalled in \cite[Theorem 2.5]{AOTree}, has
\[
 \widetilde Q_0=(Q_B)_0,\qquad
 \widetilde Q_1=(Q_B)_1\sqcup
       \{a_r:t(r)\longrightarrow s(r):r\in\mathcal R\},\qquad
 W_B=\sum_{r\in\mathcal R}r a_r,
\]
\[
 \widetilde B\iso\Jac(\widetilde Q,W_B).
\]
Every relation has length at least two, so $W_B$ has order at least three and its cyclic derivatives have order at least two. Thus the Jacobian relations remove no arrows from the radical modulo its square: the Gabriel quiver of $\widetilde B$ is $\widetilde Q$.

If $\mathcal R\ne\varnothing$, choose a nonzero path $p$ occurring with nonzero coefficient in some $r\in\mathcal R$. In the quiver $\widetilde Q$ one has
\[
 s(r)\ya{\ p\ }t(r)\ya{\ a_r\ }s(r).
\]
This is a directed closed walk, hence contains an oriented cycle. The argument concerns paths in the quiver; it does not require $p$ to survive in the Jacobian algebra. Acyclicity therefore forces $\mathcal R=\varnothing$. It follows that $I=0$, $B=\bbC Q_B$ is hereditary, and
\[
 \widetilde B=T_B\operatorname{Ext}^2_B(DB,B)=B.
\]
\end{proof}

\begin{prp}\label{prop:acyclic-lifting}
Let $A$ be $\ta_2$-finite of global dimension at most two, let $Q$ be a finite acyclic quiver, and fix a triangle equivalence
\[
 E:\calC_A\ya{\sim}\calC_{\bbC Q}.
\]
Put $U=E^{-1}(\pi_{\bbC Q}\bbC Q)$. If $U\iso\pi_A T$ for some $T\in\mathcal T_A$, then
\[
 \End_{\Db(A)}(T)\iso\bbC Q,\qquad
 \Db(A)\simeq\Db(\bbC Q).
\]
Consequently,
\[
 \left\{
 \begin{array}{l}
 A\text{ has acyclic cluster type},\\
 \pi_A(\mathcal T_A)=\CT(\calC_A)
 \end{array}
 \right.
 \quad\Longrightarrow\quad
 A\text{ is piecewise hereditary}.
\]
\end{prp}
\begin{proof}
We may take $T$ basic because $U$ is basic. Set $B=\End(T)$. Then $B$ is finite dimensional, $\operatorname{gl.dim}B\le2$, and the tilting equivalence
\[
 G:\Db(B)\ya{\sim}\Db(A),\qquad G(B)\iso T,
\]
commutes with Serre functors. Writing $\mathbb S_{2,A}=\mathbb S_A[-2]$ and similarly for $B$, we obtain from the canonical description of orbit endomorphisms that
\begin{align*}
 H^0\Pi_3(B)
 &\iso\Ds_{p\in\bbZ}\Hom_{\Db(B)}(B,\mathbb S_{2,B}^{-p}B)\\
 &\iso\Ds_{p\in\bbZ}\Hom_{\Db(A)}(T,\mathbb S_{2,A}^{-p}T)\\
 &\iso\End_{\calC_A}(\pi_A T)
 \iso\End_{\calC_{\bbC Q}}(\pi_{\bbC Q}\bbC Q)
 \iso\bbC Q.
\end{align*}
The equalities are algebra isomorphisms: the orbit products are preserved by the Serre-compatible equivalence. In particular, $B$ is $\ta_2$-finite. By \lemref{lem:acyclic-canonical}, we have $B=H^0\Pi_3(B)\iso\bbC Q$, so $G$ is a derived equivalence between $A$ and the hereditary algebra $\bbC Q$. The last assertion follows by applying the lifting hypothesis to $U$.
\end{proof}

The proof applies to every finite acyclic quiver $Q$, without a restriction on its underlying graph or representation type. In particular, if $A$ has acyclic cluster type but is not piecewise hereditary, then
\[
 E^{-1}(\pi_{\bbC Q}\bbC Q)\notin\pi_A(\mathcal T_A).
\]
The conclusion excludes every tilting lift with endomorphism algebra of global dimension at most two, not merely a chosen lift.

For comparison, \cite[Corollary 3.16]{AOTree} proves $\Db(A)\simeq\Db(\bbC Q)$ for tree cluster type without any universal-lifting assumption. For type $\widetilde A_{p,q}$, normalize $p\ge q\ge1$. Theorems 4.5 and 4.7 of \cite{AOTree} classify derived equivalence by $w$ if $p>q$ and by $|w|$ if $p=q$; Corollary 4.6 characterizes piecewise heredity by $w=0$. Hence the object above has no qualifying lift whenever $w\ne0$. For a more detailed comparison of the liftable subsets, see Proposition 6.22 of that paper. The general question follows these results in \cite{AOTree}.

\subsection{A counterexample outside acyclic cluster type}\label{subsec:torus-example}
We next give an example for which the implication fails without the acyclic-cluster-type assumption. The ordinary quiver of the algebra is acyclic.

\begin{thm}\label{thm:torus-counterexample}
Let $k=\bbC$ and let $A=kQ/(\ga a,\de c)$, where paths in this presentation are concatenated from left to right and
\[
 Q:\qquad
 4\overset{\ga,\de}{\rightrightarrows}3
   \ya{a}2\ya{b}1,
 \qquad c:3\longrightarrow1.
\]
Then $\dim_k A=13$, $\operatorname{gl.dim}A=2$, and $A$ is $\ta_2$-finite but not piecewise hereditary. Its ordinary quiver is acyclic, but $A$ is not of acyclic cluster type. Nevertheless, for every basic cluster-tilting object $U\in\calC_A$ there exists a tilting complex $T\in\Db(A)$ such that
\[
 \pi_A T\iso U,\qquad \End_{\Db(A)}(T)\iso A.
\]
In particular, $\pi_A(\mathcal T_A)=\CT(\calC_A)$.
\end{thm}

We begin with the surface realization of $A$.
The nonzero paths form the basis
\[
 e_1,e_2,e_3,e_4,\ a,b,c,\ga,\de,\ ab,\ga c,\de a,\de ab.
\]
For right projectives $P_i=e_iA$, the simple modules have projective dimensions $0,1,1,2$, respectively. More explicitly, their minimal resolutions are
\[
\begin{gathered}
 S_1=P_1,\qquad 0\longrightarrow P_1\longrightarrow P_2\longrightarrow S_2\longrightarrow0,\\
 0\longrightarrow P_2\ds P_1\longrightarrow P_3\longrightarrow S_3\longrightarrow0,\\
 0\longrightarrow P_2\ds P_1\longrightarrow P_3^{\ds2}
   \longrightarrow P_4\longrightarrow S_4\longrightarrow0.
\end{gathered}
\]
In the last sequence, $P_3^{\ds2}\to P_4$ is induced by $(\ga,\de)$; its kernel is generated by $(a,0)$ and $(0,c)$. Thus the last differential is injective and the displayed resolution is minimal.

For comparison with \cite[Section 4]{AmiotTorus}, the path algebra in that paper uses right-to-left multiplication. With that convention, $A^{\mathrm{op}}$ is its algebra $\La'_0(2,0)$: the arrow labels correspond to
\[
 a=\al_2,\qquad b=\al_1,\qquad c=\be,
\]
and the relations become $\al_2\ga$ and $\be\de$. Reversing the orientation of a marked surface reverses every arrow of its triangulation quiver and every cyclic word in its potential. Keeping the cut values therefore replaces the cut algebra by its opposite. The vanishing of the two nonseparating curve evaluations is preserved: the crossing evaluations are transported, up to the orientation signs, and zero is unchanged. By applying \cite[Proposition 4.3(1), its proof and Figure 4.1]{AmiotTorus} with $s=2$ and then reversing the surface orientation, we obtain a surface presentation
\begin{equation}
 A=A(\De_0,d_0),\qquad
 (g,b,|M|)=(1,1,1),\qquad d_0(\al)=d_0(\be)=0,
 \label{eq:torus-base}
\end{equation}
where $d_0$ is an admissible cut and $\al,\be$ are oriented simple closed curves whose classes form a basis of $H_1(S,\bbZ)$. An admissible cut assigns degree one to exactly one arrow of each internal triangle and degree zero to the remaining arrows. A degree-one map is any integer arrow grading for which each internal triangle has total degree one. The latter need not be a cut.

The surface algebra is the degree-zero part of $\Jac(Q_{\De_0},W_{\De_0},d_0)$. By \cite[Theorem 3.7]{AGSurface}, $A$ is $\ta_2$-finite and there is a triangle equivalence
\[
 \calC_A\simeq\calC_{(S,M)},\qquad \pi_A A\longmapsto U_{\De_0}.
\]
Fix this equivalence for the rest of the section. By \cite[Corollary 1.7]{BZ}, every basic cluster-tilting object is $U_\De$ for an ideal triangulation $\De$ of this same marked surface.

\subsection{Triangulations and weights on the marked torus}
We first determine the ribbon graphs associated with the triangulations.
\begin{lem}\label{lem:torus-ribbon}
For the torus with one boundary component and one marked boundary point, any two ideal triangulations are related by an orientation-preserving homeomorphism of the marked surface.
\end{lem}
\begin{proof}
Every triangulation has four internal arcs and three triangles: if $n$ and $t$ denote these numbers, Euler characteristic and incidence counting give
\[
 1-(n+1)+t=-1,\qquad 3t=2n+1,
 \qquad (n,t)=(4,3).
\]
Let $C$ denote the triangle incident with the boundary segment, and let $A,B$ denote the remaining triangles. There are no self-folded triangles: an ideal self-folded triangle encloses a puncture, and the surface has none. Consequently no dual edge is a loop. If $x_{AB},x_{AC},x_{BC}$ count internal dual edges, then
\[
 x_{AB}+x_{AC}=3,\quad x_{AB}+x_{BC}=3,\quad x_{AC}+x_{BC}=2,
\]
so $(x_{AB},x_{AC},x_{BC})=(2,1,1)$. Denote the parallel edges by $x,y$, and the edges $AC,BC$ by $u,v$. Let $l=CD$ be the edge dual to the boundary segment, where $D$ is a univalent boundary vertex. Thus $l$ determines the marked boundary point.

We determine the cyclic orders at the vertices. With half-edges indexed by their edges and incident vertices, the edge involution is
\[
 \io=(x_A\ x_B)(y_A\ y_B)(u_A\ u_C)(v_B\ v_C)(l_C\ l_D).
\]
Write the vertex rotation as
\[
 \si=(x_A\ y_A\ u_A)^{\epsilon_A}
         (x_B\ y_B\ v_B)^{\epsilon_B}
         (u_C\ v_C\ l_C)^{\epsilon_C}(l_D),
 \qquad \epsilon_A,\epsilon_B,\epsilon_C\in\{1,-1\}.
\]
The number of boundary components of the thickened ribbon graph is the number of cycles of $\si\io$. By composing the permutations, we obtain
\[
\begin{array}{c|cc}
 (\epsilon_A,\epsilon_B)&\epsilon_C=1&\epsilon_C=-1\\\hline
 (1,1)&1&1\\
 (1,-1)&3&3\\
 (-1,1)&3&3\\
 (-1,-1)&1&1
\end{array}
\]
For example, the $(1,1,1)$ case has the single boundary cycle
\[
 (x_A\ y_B\ u_A\ v_C\ x_B\ y_A\ v_B\ l_C\ l_D\ u_C).
\]
Since the Euler characteristic is $4-5=-1$, the genus $g$ and
the number $b$ of boundary components satisfy $2g+b=3$. Hence
\[
 (g,b)=
 \begin{cases}
 (1,1),&\epsilon_A=\epsilon_B,\\
 (0,3),&\epsilon_A\ne\epsilon_B.
 \end{cases}
\]

Put $\epsilon=(\epsilon_A,\epsilon_B,\epsilon_C)$ and write
$\si_\epsilon$ for the corresponding vertex rotation. Consider the
following permutations of the half-edges:
\[
\begin{aligned}
 \tau_1&=(x_A\ y_A)(x_B\ y_B),\\
 \tau_2&=(x_A\ x_B)(y_A\ y_B)(u_A\ v_B)(u_C\ v_C).
\end{aligned}
\]
For $i=1,2$, we have
\[
 \tau_i\io\tau_i^{-1}=\io,\qquad
 \tau_i(l_C)=l_C,\qquad \tau_i(l_D)=l_D.
\]
Moreover,
\[
\begin{aligned}
 \tau_1\si_\epsilon\tau_1^{-1}
   &=\si_{(-\epsilon_A,-\epsilon_B,\epsilon_C)},\\
 \tau_2\si_\epsilon\tau_2^{-1}
   &=\si_{(\epsilon_B,\epsilon_A,-\epsilon_C)}.
\end{aligned}
\]
Thus $\langle\tau_1,\tau_2\rangle$ acts transitively on
\[
 \{\si_{(a,a,c)}:a,c\in\{1,-1\}\}.
\]
The four ribbon graphs are therefore isomorphic by isomorphisms fixing
the boundary edge. Since these isomorphisms preserve the cyclic orders,
they extend to orientation-preserving homeomorphisms of their
thickenings. The extensions identify the corresponding triangles and
paired sides, and hence the triangulations. They also preserve the
marked boundary point, as required.
\end{proof}

We next compare the weights of degree-one maps on the fixed marked surface.
Fix the curves $\al,\be$ in \eqref{eq:torus-base} and the positively oriented boundary curve $c_\partial$. Write
\[
 w(d)=\bigl(d(\al),d(\be),d(c_\partial)\bigr).
\]
The evaluations are the signed crossing evaluations of \cite[Section 2.2]{AGSurface}; the curves remain fixed when the triangulation changes. For every degree-one map on this surface,
\[
 d(c_\partial)=4g-4+2b=2
 \qquad\text{by \cite[Proposition 2.9]{AGSurface}}.
\]
In particular, $d$ itself does not define a homomorphism on all of $H_1(S,\bbZ)$: the boundary class is zero but its evaluation is two.

\begin{lem}\label{lem:torus-marked-weights}
Let $\De$ be any ideal triangulation. Choose an orientation-preserving homeomorphism $\Ph$ with $\Ph\De_0=\De$, and transport the cut by
\[
 d_\De(\Ph a)=d_0(a).
\]
Let $\de_\De$ be obtained from $d_0$ by left graded mutations along any flip sequence from $\De_0$ to $\De$. Then
\[
 w(d_\De)=w(\de_\De)=(0,0,2),\qquad
 \de_\De-d_\De=\partial_0^*r
 \quad\text{for some }r\in\bbZ^{(Q_\De)_0}.
\]
Here $(\partial_0^*r)(a)=r(t(a))-r(s(a))$.
\end{lem}
\begin{proof}
By naturality of the crossing evaluation \cite[Lemma 2.6]{AGSurface}, we have
\[
 d_\De(\et)=d_0(\Ph^{-1}\et).
\]
The homeomorphism $\Ph$ preserves simplicity and induces an automorphism of $H_1(S,\bbZ)$. Thus both $\Ph^{-1}\al$ and $\Ph^{-1}\be$ are simple closed curves with nonzero homology class. By \cite[Proposition 3.11]{AmiotTorus}, for every such curve $\et$ we have
\[
 [\et]=m[\al]+n[\be]\quad\Longrightarrow\quad
 d_0(\et)=m d_0(\al)+n d_0(\be)=0.
\]
It follows that $w(d_\De)=(0,0,2)$. Since $\Ph$ identifies the internal triangles and their chosen cut arrows, $d_\De$ is an admissible cut.

On the other hand, \cite[Lemma 2.14]{AGSurface} preserves the weights of the fixed curves under each graded flip. By induction along the chosen sequence, we obtain
\[
 w(\de_\De)=w(d_0)=(0,0,2).
\]
Both $\de_\De$ and $d_\De$ are degree-one maps. Thus \cite[Corollary 2.11]{AGSurface} applies. To see the resulting coboundary explicitly, consider the complex
\[
 \bbZ^{(Q_\De)_0}\ya{\partial_0^*}
 \bbZ^{(Q_\De)_1}\ya{\partial_1^*}
 \bbZ^{(Q_\De)_2};
\]
Set $z=\de_\De-d_\De$. The degree-one conditions and the two nonseparating periods give
\[
 \partial_1^*z=0,\qquad z(\al)=z(\be)=0.
\]
By \cite[Remark 2.7]{AGSurface}, evaluation on these two periods induces
\[
 H^1(C^\bullet)\ya{\sim}\bbZ^2,
 \qquad [z]\longmapsto\bigl(z(\al),z(\be)\bigr).
\]
Hence $[z]=0$ and $z=\partial_0^*r$. Here the cocycle is the difference of the two degree-one maps, and the evaluations are taken on the same fixed curves $\al,\be$. Equality up to a change of homology basis or equality of the gcd invariants would not suffice.
\end{proof}

\subsection{Tilting lifts}
\begin{prp}\label{prop:torus-marked-lift}
Every basic $U_\De\in\calC_A$ has a tilting lift $T_\De$ with
\[
 \pi_A T_\De\iso U_\De,\qquad
 \End(T_\De)\iso A(\De,d_\De)\iso A.
\]
\end{prp}
\begin{proof}
Put $F_A=\mathbb S_A[-2]$. Take a finite flip sequence from $\De_0$ to $\De$ and perform the corresponding left mutations in $\Db(A)$, as in \cite[Definition 3.8]{AGSurface}. Write the resulting object as $X=\Ds_iX_i$, with summands labelled by the arcs of $\De$. The mutation triangles and \cite[Theorem 3.9]{AGSurface} give
\[
 \pi_A X=U_\De,\qquad
 \Ga_X:=\Ds_{p\in\bbZ}\Hom_{\Db(A)}(X,F_A^{-p}X)
 \iso\Jac(Q_\De,W_\De,\de_\De)
\]
as graded algebras with the arc idempotents retained. At this point $X$ is not asserted to be tilting.

Set $B=A(\De,d_\De)$. By \cite[Theorem 3.7]{AGSurface}, we have
\[
\begin{gathered}
 \operatorname{gl.dim}B\le2,\qquad B\text{ is }\ta_2\text{-finite},\\
 \Ga_B:=\Ds_p\Hom_{\Db(B)}(B,F_B^{-p}B)
 \iso\Jac(Q_\De,W_\De,d_\De).
\end{gathered}
\]
By \lemref{lem:torus-marked-weights}, $\de_\De-d_\De=\partial_0^*r$.
For a path $a_\ell\cdots a_1:i\to j$, the degree difference is
\[
 \sum_{\nu=1}^{\ell}\bigl(r(t(a_\nu))-r(s(a_\nu))\bigr)
 =r(j)-r(i).
\]
Vertex shifts therefore identify the covering categories, preserving
composition, relations and the translation action.

Let
\[
 \mathcal V_X=\add\{F_A^pX_i:i\in(Q_\De)_0,\ p\in\bbZ\}
\]
be the full additive subcategory of $\Db(A)$ generated by these objects. It is the covering category of $\Ga_X$, since
\[
 \Hom(F_A^pX_i,F_A^qX_j)
   \iso\Hom(X_i,F_A^{q-p}X_j).
\]
Under the equivalence of covering categories, the projective $e_iB$
corresponds to $T_i=F_A^{h_i}X_i$ for an integer $h_i$. Put $T=\Ds_iT_i$.
With the convention that an arrow $a:i\to j$ of degree $d(a)$ represents
a morphism $X_i\to F_A^{-d(a)}X_j$, we have
\[
 F_A^{h_i}X_i\longrightarrow
 F_A^{-d(a)+h_i-h_j}(F_A^{h_j}X_j),
 \qquad d'(a)=d(a)+h_j-h_i.
\]
Since $\de_\De-d_\De=\partial_0^*r$, the choice $h_i=-r(i)$
identifies the two gradings. The vertex labels are unchanged.

Consequently, there is an additive equivalence commuting with $F$
\[
 f:\mathcal U_B:=\add\{F_B^pB:p\in\bbZ\}\ya{\sim}\mathcal V_X,
 \qquad f(B)=T,\qquad \pi_A T=U_\De.
\]
We claim that
\[
 \mathcal V_X=\pi_A^{-1}(\add U_\De).
\]
Since $\pi_A F_A\simeq\pi_A$ and $\pi_A X=U_\De$, the inclusion from left to right follows. Conversely, let $\pi_A Y$ be a summand of $\pi_A X^{\ds m}$. Choose split maps $\bar a:\pi_A Y\to\pi_A X^{\ds m}$ and $\bar b:\pi_A X^{\ds m}\to\pi_A Y$ with $\bar b\bar a=1$. Full faithfulness of the orbit category in its triangulated hull writes these as finite sums
\[
 a_p:Y\longrightarrow F_A^pX^{\ds m},\qquad
 b_q:X^{\ds m}\longrightarrow F_A^qY.
\]
The orbit-degree-zero component of $\bar b\bar a=1$ is
\[
 \sum_p F_A^p(b_{-p})a_p=1_Y.
\]
Thus the maps $(a_p)_p$ and $\bigl(F_A^p(b_{-p})\bigr)_p$ exhibit $Y$ as a retract of the finite sum $\Ds_pF_A^pX^{\ds m}$. Hence $Y\in\mathcal V_X$. By \cite[Proposition 3.1]{AO}, this inverse image is a cluster-tilting subcategory of $\Db(A)$.

Since $\Db(A)$ is algebraic with Serre functor, $\mathcal V_X$ is cluster-tilting, and $f$ commutes with $F=\mathbb S[-2]$, \cite[Theorem 3.5]{AO} applies. Hence there exists a triangle equivalence
\[
 G:\Db(B)\ya{\sim}\Db(A),\qquad G|_{\mathcal U_B}\iso f.
\]
In particular, $G(B)\iso T$. Therefore
\[
\begin{gathered}
 \Hom(T,T[n])\iso\Hom(B,B[n])=0\quad(n\ne0),\\
 \thick(T)=\Db(A),\qquad \End(T)\iso B.
\end{gathered}
\]
Finally, $\Ph$ identifies the cut quivers with relations, so $B\iso A$. By construction, $\pi_A T=U_\De$.
\end{proof}

\subsection{Failure of piecewise heredity}
\begin{lem}\label{lem:torus-not-ph}
The algebra $A$ in \thmref{thm:torus-counterexample} is not piecewise hereditary.
\end{lem}
\begin{proof}
The path-count matrix, with row $i$ recording paths starting at $i$, and its inverse are
\[
 C=\begin{pmatrix}1&0&0&0\\1&1&0&0\\2&1&1&0\\2&1&2&1\end{pmatrix},
 \qquad E=C^{-1}=\begin{pmatrix}1&0&0&0\\-1&1&0&0\\-1&-1&1&0\\1&1&-2&1\end{pmatrix}.
\]
Transposing the Cartan convention does not change the following invariants. By direct calculation, we obtain
\[
 \det(xI+C^{-t}C)=(x-1)^2(x^2+x+1),\qquad
 \operatorname{tr}(-C^{-t}C)=1,
\]
\[
 E+E^t=\begin{pmatrix}2&-1&-1&1\\-1&2&-1&1\\-1&-1&2&-2\\1&1&-2&2\end{pmatrix},
 \qquad \operatorname{Spec}(E+E^t)=\{5,3,0,0\}.
\]
Thus the symmetrized Euler form is positive semidefinite of corank two.

Suppose $A$ were piecewise hereditary. Since it is connected and its Coxeter trace is greater than $-1$, Happel's criterion, as stated in \cite[Theorem 1.1(1)]{LRT}, implies $\Db(A)\simeq\Db(kR)$ for a connected acyclic quiver $R$. A derived equivalence induces an integral isometry of Euler forms. Consequently $2I-N_R$ must be positive semidefinite of corank two, where $N_R$ is the symmetric adjacency matrix counting arrows in either direction. Positive semidefiniteness implies that every eigenvalue of $N_R$ is at most two; singularity implies that two is an eigenvalue. Since $N_R$ is a connected nonnegative symmetric matrix, Perron--Frobenius makes its largest eigenvalue simple. Hence $2I-N_R$ has corank one, a contradiction.
\end{proof}

\begin{proof}[Proof of \thmref{thm:torus-counterexample}]
By the basis and projective resolutions above, we obtain the stated dimension and global dimension. By the surface realization, the algebra is $\ta_2$-finite. \prpref{prop:torus-marked-lift} treats every basic cluster-tilting object, while \lemref{lem:torus-not-ph} excludes piecewise heredity. If $A$ had acyclic cluster type, \prpref{prop:acyclic-lifting} would then force it to be piecewise hereditary, a contradiction. Its displayed ordinary quiver is nevertheless acyclic. The reverse inclusion $\pi_A(\mathcal T_A)\subseteq\CT(\calC_A)$ is the usual tilting-to-cluster-tilting construction; alternatively it follows by applying the derived equivalence associated with a tilting complex whose endomorphism algebra has global dimension at most two. This proves the equality in the statement.
\end{proof}

\begin{rmk}
The surface argument uses two properties of the marked torus: the transitivity of orientation-preserving homeomorphisms on its triangulations, and the vanishing of the chosen cut on nonseparating simple closed curves. For another surface, both properties require separate verification.
\end{rmk}


\begin{thebibliography}{99}
\bibitem{AI}T. Aihara and O. Iyama, \emph{Silting mutation in triangulated categories}, Journal of the London Mathematical Society (2) \textbf{85} (2012), no.~3, 633--668.

\bibitem{Amiot}C. Amiot, \emph{Cluster categories for algebras of global dimension 2 and quivers with potential}, Annales de l'Institut Fourier \textbf{59} (2009), no.~6, 2525--2590. 

\bibitem{AmiotSurvey}C. Amiot, \emph{On generalized cluster categories}, in \emph{Representations of Algebras and Related Topics}, EMS Series of Congress Reports, European Mathematical Society, 2011, 1--53. 


\bibitem{AmiotTorus}C. Amiot, \emph{The derived category of surface algebras: the case of the torus with one boundary component}, preprint, version 3 (25 March 2016). \href{https://arxiv.org/abs/1506.02410v3}{arXiv:1506.02410v3}. 

\bibitem{AGSurface}C. Amiot and Y. Grimeland, \emph{Derived invariants for surface algebras}, preprint, version 2 (25 March 2016). \href{https://arxiv.org/abs/1411.0383v2}{arXiv:1411.0383v2}. 

\bibitem{AOTree}C. Amiot and S. Oppermann, \emph{Algebras of acyclic cluster type: tree type and type $\widetilde A$}, Nagoya Mathematical Journal \textbf{211} (2013), 1--50. 

\bibitem{AO}C. Amiot and S. Oppermann, \emph{Cluster equivalence and graded derived equivalence}, Documenta Mathematica \textbf{19} (2014), 1155--1206.

\bibitem{AP}H. Asashiba and S. Pan, \emph{Cohen--Montgomery duality for bimodules and singular equivalences of Morita type}, preprint, version 4 (2026). \href{https://arxiv.org/abs/2408.03280v4}{arXiv:2408.03280v4}.

\bibitem{BZ}T. Br\"ustle and J. Zhang, \emph{On the cluster category of a marked surface without punctures}, Algebra \& Number Theory \textbf{5} (2011), no.~4, 529--566.

\bibitem{COS}A. Canonaco, M. Ornaghi and P. Stellari, \emph{Localizations of the category of $A_\infty$ categories and internal Homs}, Documenta Mathematica \textbf{24} (2019), 2463--2492. 

\bibitem{SGA}M. Demazure and A. Grothendieck (eds.), \emph{Sch\'emas en groupes (SGA 3)}, Tomes I--II, Lecture Notes in Mathematics \textbf{151--152}, Springer, 1970.

\bibitem{DWZ}H. Derksen, J. Weyman and A. Zelevinsky, \emph{Quivers with potentials and their representations I: Mutations}, Selecta Mathematica (N.S.) \textbf{14} (2008), 59--119.

\bibitem{Fushimi}R. Fushimi, \emph{Contravariant Koszul duality between non-positive and positive dg algebras}, preprint, version 2 (17 February 2025). \href{https://arxiv.org/abs/2409.08842v2}{arXiv:2409.08842v2}. 

\bibitem{Ginzburg}V. Ginzburg, \emph{Calabi--Yau algebras}, preprint (2006), version 3 (9 January 2007). \href{https://arxiv.org/abs/math/0612139v3}{arXiv:math/0612139v3}. 

\bibitem{Han}N. Hanihara, \emph{Calabi--Yau completions for roots of dualizing dg bimodules}, preprint, version 1 (2024). \href{https://arxiv.org/abs/2412.18753v1}{arXiv:2412.18753v1}.

\bibitem{HI}N. Hanihara and O. Iyama, \emph{Silting correspondences and Calabi--Yau dg algebras}, preprint, version 2 (2025). \href{https://arxiv.org/abs/2508.12836v2}{arXiv:2508.12836v2}. 

\bibitem{HK}Z. Hua and B. Keller, \emph{Cluster categories and rational curves}, Geometry \& Topology \textbf{28} (2024), 2569--2634. 

\bibitem{HZ}Z. Hua and G.-S. Zhou, \emph{Noncommutative Mather--Yau theorem and its applications to Calabi--Yau algebras}, Mathematische Annalen \textbf{386} (2023), 1605--1639. 

\bibitem{HZFinite}Z. Hua and G. Zhou, \emph{Finite determinacy and rigidity of quiver potentials}, preprint, version 1 (26 October 2026). \href{https://arxiv.org/abs/2610.35466v1}{arXiv:2610.35466v1}. 

\bibitem{HZS}B. Huisgen-Zimmermann and M. Saor\'in, \emph{Geometry of chain complexes and outer automorphisms under derived equivalence}, Transactions of the American Mathematical Society \textbf{353} (2001), no.~12, 4757--4777. 


\bibitem{KalckYang}M. Kalck and D. Yang, \emph{Derived categories of graded gentle one-cycle algebras}, Journal of Pure and Applied Algebra \textbf{222} (2018), no.~10, 3005--3035. 

\bibitem{KellerAinf}B. Keller, \emph{Introduction to $A$-infinity algebras and modules}, Homology, Homotopy and Applications \textbf{3} (2001), no.~1, 1--35.

\bibitem{KellerOrbit}B. Keller, \emph{On triangulated orbit categories}, Documenta Mathematica \textbf{10} (2005), 551--581. 

\bibitem{KellerDG}B. Keller, \emph{On differential graded categories}, in \emph{International Congress of Mathematicians}, Vol.~II, European Mathematical Society, 2006, 151--190.

\bibitem{KellerOrbitCorr}B. Keller, \emph{Corrections to ``On triangulated orbit categories''}, author's note, revised 21 September 2009. \href{https://webusers.imj-prg.fr/~bernhard.keller/publ/corrTriaOrbit.pdf}. 

\bibitem{Keller}B. Keller, with an appendix by M. Van den Bergh, \emph{Deformed Calabi--Yau completions}, Journal f\"ur die reine und angewandte Mathematik \textbf{654} (2011), 125--180. 

\bibitem{KY}B. Keller and D. Yang, \emph{Derived equivalences from mutations of quivers with potential}, Advances in Mathematics \textbf{226} (2011), no.~3, 2118--2168.

\bibitem{LRT}M. Lanzilotta, M. J. Redondo and R. Taillefer, \emph{Combinatorial classification of piecewise hereditary algebras}, preprint, version 1 (25 March 2009). \href{https://arxiv.org/abs/0903.4326v1}{arXiv:0903.4326v1}. 

\bibitem{Rickard}J. Rickard, \emph{Morita theory for derived categories}, Journal of the London Mathematical Society (2) \textbf{39} (1989), no.~3, 436--456. 

\bibitem{VdB}M. Van den Bergh, \emph{Calabi--Yau algebras and superpotentials}, Selecta Mathematica (N.S.) \textbf{21} (2015), 555--603. 

\bibitem{XYZ}Z. Xie, D. Yang and H. Zhang, \emph{A differential graded approach to the silting theorem}, Journal of Pure and Applied Algebra \textbf{227} (2023), no.~2, 107180. 

\end{thebibliography}
\end{document}